\documentclass[12pt]{amsart} %
\usepackage{verbatim}
\usepackage{multirow}
\usepackage[bottom]{footmisc}
\usepackage{graphicx}

\def\rat#1{\mathcal R_{#1}}

\def\Forall{\quad \hbox{ for all }}
\def\bal#1\eal{\begin{aligned} #1 \end{aligned}}
\def\beq#1\eeq{\begin{equation} #1 \end{equation}}

\def\argmin{\mathop{\rm argmin}}
\def\bfu{{\bf u}}

\def\bff{{\bf f}}
\def\calA{{\mathcal A}}

\def\calT{{\mathcal T}}

\def\calAt{{\mathcal A}_h}
\def\calBt{{\mathcal B}_h}

\def\calIt{{\mathcal I}_h}
\def\calTt{{\mathcal T}_h}
\def\D{{M}}

\def\wcalSt{{\mathbb S}}
\def\wcalMt{{\mathbb M}}
\def\wcalAt{ {\mathbb A}}

\def\wcalIt{ {\mathbb I}}

\def\ttt{\tilde t}
\def\tkk{\tilde k}
\def\D{N}

\def\bfpsi{{\boldsymbol \psi}}

\newcommand{\RR}{{\mathbb{R}}}

\def\argmin{\mathop{\rm argmin}}

\def\tiluh{{ \bf u}}
\def\tilfh{{ \bf f}}
\def\tilwh{{\bf u}_k}
\def\tilvh{{\bf v}}
\def\tilF{{\bf F}}

\def\calI{\mathcal I}

\def\tt{\xi}
\def\q{c}
\def\qq{q}
\def\kk{{k^\prime}}

\def\zz{z}
\def\trg{\widetilde r_{\alpha,k}}
\def\Uhxt{{U}_h}  % this replaces  $\tilde{u}_{h}$

\newcommand{\R}{\mathbb{R}}

\newcommand{\II}{\mathbb{I}}

\newcommand{\rqa}{r_{q,\alpha,k}}
\newcommand{\brqa}{{r}_{q_1,q_2,\alpha,k}}
\newcommand{\Eqa}{E_{q,\alpha,k}}
\newcommand{\bEqa}{\bar{E}_{q_1,q_2,\alpha,k}}

\def\argmin{\mathop{\rm argmin}}

\numberwithin{equation}{section}
\numberwithin{figure}{section}
\theoremstyle{plain}% default
\newtheorem{theorem}{Theorem}[section]
\newtheorem{lemma}[theorem]{Lemma}

\newtheorem{remark}[theorem]{Remark}

 \usepackage[pdfpagelabels]{hyperref} 
 \title[A Survey: Numerics for Spectral Elliptic PDEs]{A Survey on Numerical Methods for  Spectral
 Space-Fractional Diffusion Problems}
 \author%[S. Harizanov, R. Lazarov, S. Margenov]
 {Stanislav Harizanov}
 \address{Institute of Information and Communication Technologies\\
Bulgarian Academy of Sciences\\
Acad. G. Bontchev Str., Block 25A\\
1113 Sofia, BULGARIA}
\email{sharizanov@parallel.bas.bg}
 
 \author{ Raytcho Lazarov}
 \address{Texas A\&M University\\
Department of Mathematics \\
505 D Blocker Build.\\
77843 College Station, TX, USA}
\email{lazarov@math.tamu.edu}
 
 \author{ Svetozar Margenov}
 \address{Institute of Information and Communication Technologies\\
Bulgarian Academy of Sciences\\
Acad. G. Bontchev Str., Block 25A\\
1113 Sofia, BULGARIA}
\email{margenov@parallel.bas.bg}

\begin{document}

% \vbox to 2.0cm { \vfill }

%%% to make empty space of approx. 2 cm %%%%%%
%%% will be replaced by Editor with the journal's and publoishers logos %%%%%%%%

 \bigskip \medskip

%%%% Abstract %%%%%%%%%%%%%%%%%%%%%%%%%
 \begin{abstract}
The survey is devoted to numerical solution of the equation $\calA^\alpha u=f$, $0 < \alpha <1$,
where $\calA$ is a symmetric positive definite operator corresponding to a second order 
elliptic boundary value problem in a bounded domain $\Omega$ in $\R^d$.
The fractional power  $\calA^\alpha $ is a non-local operator and is defined though the spectrum of  $\calA$.
Due to growing interest and demand in applications of sub-diffusion models to physics and engineering, 
in the last decade, several numerical approaches have been proposed, studied, and tested. 
We consider discretizations of the elliptic operator $ \calA $ by using 
an $N$-dimensional finite element space $V_h$ or finite differences over a uniform mesh with $N$ grid points.  
In the case of finite element  approximation we get a 
symmetric and positive definite operator $\calAt: V_h \to V_h$, which
results in an operator equation $\calAt^{\alpha} u_h = f_h$ for $u_h \in V_h$.

The numerical solution of this equation is based on the following three equivalent representations of the solution: 
(1) Dunford-Taylor integral formula (or its equivalent Balakrishnan formula, \eqref{bal}), 
(2)  extension of the a second order elliptic problem in $\Omega \times (0,\infty)\subset \R^{d+1}$ 
\cite{caffarelli2007extend,nochetto2015pde} 
(with a local operator) or as a pseudo-parabolic equation in the
cylinder $(x,t) \in \Omega \times (0,1) $,
\cite{Vabishchevich14, DLP-Pade},
(3) spectral representation \eqref{eq:spectral-0} and the best uniform rational approximation (BURA) of $z^\alpha$
on $[0,1]$, \cite{HLMMV18,harizanov2020JCP}.
Though substantially different in  origin and their analysis, these methods can be interpreted as some rational approximation
of $\calAt^{-\alpha}$. In this paper we present the main ideas of these methods and the corresponding algorithms, 
discuss their accuracy, computational complexity and compare their efficiency and robustness.

\medskip

\smallskip


\end{abstract}

\maketitle

%%%%%%% end make title %%%%%%%%%%%%%%%%%%%%%%%%%%%%%%%%%%
% \vspace*{-20pt}

%%%%%%%% begin papers' body %%%%%%%%%%%%%%%%%%%%%%%%%%%%%

%%%%%%%%%%%%%%%%%%%%%%%%%%% Section 1 %%%%%%%%%%%%%

\section{Introduction}\label{sec:1}
\setcounter{section}{1}
\setcounter{equation}{0}\setcounter{theorem}{0}

Fractional calculus is an emerging field in mathematics.  
Equations involving fractional partial derivatives are systematically used  
to model anomalous processes
in which the Brownian motion hypotheses are violated.
The rapidly increasing interest in development of
efficient numerical methods for such problems is
motivated by the great capacity of such mathematical models in 
applications of anomalous diffusion to  
science and engineering. A collection of such 
real world applications is presented in \cite{SZBCC-18} where
experts in various fields of science and engineering 
presented applied problems in Physics, Control, Signal and Image Processing, 
Mechanics and Dynamic Systems, Biology, 
Environmental Science, Materials, Economic, and
Multidisciplinary in Engineering Fields. 

In mathematics and physics, fractional order differential operators 
appear naturally in trace theory of functions in Sobolev classes 
(Sobolev embedding) \cite{BKM-19}, the theory of special classes of analytic functions, \cite{Djrbashian:1993}, 
Caputo and Riemann-Liouville fractional derivatives, 
\cite{Podlubny:1999}.  The importance of the field is demonstrated by its capabilities in 
modeling various real life phenomena, e.g. particle movement in heterogeneous media, 
\cite{MetzlerKlafter:2004} and/or  heavily tailed Levy flights of particles, \cite{HatanoHatano:1998}, 
peridynamics (deformable media with fractures), \cite{peridynamics}, image reconstruction, 
\cite{gilboa2008nonlocal},  transport of CO2 in heterogeneous media, \cite{CH-18},
phase-field crystal modeling, \cite{AM-20},
etc. In the applications there are substantial variations that involve both, 
transient and steady-state problems. For example, there are models
with fractional time derivatives of Caputo or Riemann-Liouville, 
\cite{Podlubny:1999}, and steady-state sub-diffusion problems involving fractional Laplacian, 
\cite{LPGSGZMCMAK-18}. 
The most important property of these operators is that they are non-local.

This survey is devoted to numerical methods for solving the problem $\calA^\alpha u =f$,  
where $\calA$ is an elliptic operator of second order and $0< \alpha <1$.
The simplest example of such a problem is 
the {\it spectral fractional Laplacian} (more general elliptic operators are discussed in Section \ref{ss:problem}), 
defined through the spectrum $(\lambda_j, \psi_j(x))$ of the Laplace operator $-\Delta=\calA$ 
\beq\label{fr-Laplace}
(-\Delta)^\alpha u (x)= \sum_{j=1}^\infty \lambda_j^{\alpha} (u,\psi_j) \psi_j (x) , \ \forall \ x \in \Omega,
\eeq
for functions that satisfy $u(x)=0$ on $\partial \Omega$ and
$  
 \sum_{j=1}^\infty \lambda_j^{2\alpha} |(u,\psi_j)|^2 < \infty .
$
Then the corresponding boundary value problem is: for $f \in L^2(\Omega) $ find $u$ such that
\beq\label{int-Laplace}
(-\Delta )^\alpha u (x) =  f(x) \ \ \mbox{in} \  \  \Omega,  \ \ u=0 \ \ \mbox{on} \ \ \partial \Omega .
\eeq

The {\it integral fractional Laplacian} represents another 
class of nonlocal operators. In strong sense, it is introduced via the relation:
\beq\label{integral}
(-\Delta )^\alpha u (x) = C_{d,\alpha} P.V. \int_{R^d} \frac{ u(x) -u(y) }{|x - y |^{d+2\alpha}} dy,  \ \forall \ x \in \RR^d,
\eeq
where $C_{d,\alpha}= 2^{2\alpha} \alpha \Gamma(\alpha +d/2)/(\pi^{d/2} \Gamma(1-\alpha))$ 
(e.g. \cite[formula (3.1)]{BBNOS-18,DU2020numerical}). 
Here  $ (-\Delta )^\alpha$ acts on  the set of functions that are extended by zero 
to $\RR^d$. Thus, the corresponding boundary value problem is: find $u$ s.t.
$$
(-\Delta )^\alpha u (x) =  f(x) \ \ \mbox{in} \  \  \Omega,  \ \ u=0 \ \ \mbox{in} \ \Omega^c = \RR^d \setminus \Omega  .
$$
For the corresponding weak formulation we refer to \cite{bonito2019numerical}. For the related evolution 
problem we could refer, e.g. to \cite{ABB-19}.
This problem has a probabilistic interpretation in particles random walk with arbitrary long jumps. 

There are also other definitions of fractional Laplacian that include Balakrishnan formula \eqref{bal},
formula involving semi-group, Dynkin's definition based on probabilistic considerations, e.g. \cite{DU2020numerical,kwasnicki2017ten,LPGSGZMCMAK-18}.  
As shown in  \cite{kwasnicki2017ten}
these are all equivalent in the whole space $\RR^d$ and differ substantially 
when considered in a bounded domain. Here we shall follow the spectral definition, 
discussed in details in Sections \ref{sec:2}.

The discretization of the problem $\calA^{\alpha} u=f$ is done via approximation of 
the differential operator $\calA$ by finite elements or finite differences resulting in a 
symmetric matrix $\wcalAt$ that acts on the vector $\tiluh$ of the 
unknown values  of $u$ at the mesh points (for more details, see Section \ref{sec:3}).
Then the desired approximation is $\wcalAt^\alpha  \tiluh= \tilfh$.
Here $\wcalAt^\alpha$ is a particular case of 
the general definition of a
function $f(\wcalAt)$ of the matrix $\wcalAt$, e.g. \cite{higham2008functions}, given by 
the Cauchy integral formula 
\begin{equation}\label{matrixfunction}
f(\wcalAt) = \frac{1}{2\pi i}\int_\Gamma f(\mu) (\mu \wcalIt - \wcalAt)^{-1} d \mu,
\end{equation}
where $\Gamma$ is a closed contour  that lies in the region of analyticity of $f$  and winds 
once around the spectrum in the anti-clockwise direction. In the particular case of symmetric matrix   
$\wcalAt$  and $f(\wcalAt)= \wcalAt^\alpha$ we get \eqref{bal}.

In general, this representation is not always useful from a computational point of view,  as it requires information 
about the spectral region of $\wcalAt$. Nevertheless, it is a good starting point for developing various numerical methods 
for approximate computing of $f(\wcalAt)$. For diagonalizable matrices one can use the Symbolic Math Toolbox of
Matlab R2008b, \cite{mathworks09}.  There one can find a number of algorithms developed for computing 
$p$ root, exponent, logarithm, etc of square matrices, see, e.g. \cite{higham2008functions,mathworks09}. 
However, all these methods 
are efficient for matrices of up to a moderate size.

In the case of $\wcalAt$ being a symmetric and positive definite matrix, 
the above formula has various simplified forms. 
One of them, formula \eqref{bal}, has been used to derive efficient algorithms for computing 
$\wcalAt^{-\alpha} {\bf f}$  when the corresponding matrix $\wcalAt$ is sparse and of very large size,
e.g. \cite{bonito2019numerical,bonito2019sinc,BP-15,BP-17}.

Aimed at presumably more realistic applications, we are 
interested in numerical methods for spectral fractional diffusion 
problems in multidimensional domains with  general geometry,
which after proper discretization produce large sparse symmetric matrices.  This in 
particular means that methods based on Fast Fourier Transform (FFT) 
for problems with constant coefficients in domains that are tensor product of intervals 
are outside the scope of this paper.  

The direct application of the spectral decomposition of $\wcalAt^\alpha$ involves 
computation of the eigenvalues and eigenfunctions of $\wcalAt$.
Generally, this is unacceptably expensive in terms of computations and 
computer memory requirement.  Nevertheless, such approach 
could be made quite efficient in the case of approximation of the 
elliptic operator by a spectral numerical method, supposing that the 
size of $\wcalAt$ is small enough, \cite{SXK-17}. A key point to the 
achieved effectiveness in this paper is the assumption for high 
smoothness of the solution. However, in this survey we target 
a much more general class of problems in complex domains where 
the solution could be of a very low regularity. All these  
naturally lead to large-scale linear systems with number of
unknowns in the range of hundreds of thousands and hundreds of
millions.
The methods that will be discussed from now on avoid  the explicit 
use of $\wcalAt^\alpha$  including matrix vector multiplication.

We start our discussion with the basic problem in linear algebra, namely the derivation of
solvers for linear systems with dense matrices.
In \cite{druskin1998extended}
an  extended Krylov subspace method based on the subspace 
$$K^{k,m}(\wcalAt,\phi) = 
span \{\wcalAt^{-k+1}\phi,
\hdots, \wcalAt^{-1}\phi, \phi, \hdots, \wcalAt^{m-1}\phi \}, \ m\ge 1, k\ge 1.
$$
has been proposed. As an alternative, in \cite{ilic2009numerical}, an adaptively preconditioned thick 
restart Lanczos  procedure is applied to the system with $\wcalAt$ first.  
The gathered spectral information is then used to solve the system 
with $\wcalAt^\alpha$. 
Both methods have shown significant progress.
However, they are not robust with respect to the condition number of
$\wcalAt$ and show substantial increase in the needed computer memory for ill-conditioned matrices.

In the last decade a number of new  
approaches for numerical solution of non-local fractional diffusion problems were proposed,
justified and tested. Among these are methods based on:
\begin{enumerate}
\item extension of the problem from $\Omega\subset \R^d$ to an
elliptic problem in $\Omega \times (0,\infty)\subset \R^{d+1}$ 
\cite{caffarelli2007extend,nochetto2015pde,nochetto2016pde} (with a local operator 
or a reformulation as a pseudo-parabolic problem in the
cylinder $(x,t) \in \Omega \times (0,1) $,
\cite{Vabishchevich14,Vabishchevich15,CV-20,DLP-Pade,Vabishchevich18}; 
\item  methods based on approximation of the Dunford-Taylor 
integral representation of the solution, \cite{bonito2019numerical,bonito2019sinc,BP-15,BP-17};
\item methods based on the best uniform rational approximation (BURA) of $z^\alpha$
on $[0,1]$, \cite{HLMMV18,harizanov2019cmwa,harizanov2020JCP}. 
\end{enumerate}
The scope of this survey is the original formulation of these methods
and their further development and extensions. We stress, that these methods are 
derived in different ways and employ different theoretical analysis.
However, they are interrelated and they
all can be interpreted as rational approximations of $\calAt^{-\alpha}$
or $\wcalAt^{-\alpha}$, see, e.g. \cite{
HOFREITHER2019}, that provides a solid basis for comparison and evaluation.  
Among the discussed below  properties are the exponential
convergence with respect to the degree of rational approximation, the
robustness with respect to the condition number $\kappa(\wcalAt)$,
and the nearly optimal computational complexity $O(N\log N)$, where
$N$ is the number of the unknowns in the discrete problem.

The paper is organized as follows. 
The spectral space-fractional diffusion problems are defined 
in Section~\ref{sec:2}, including the spectral fractional powers 
of an elliptic operator  $\calA$, the sub-diffusion-reaction problems and the
basic regularity properties.
The finite element and finite difference discretizations of 
$\calA$ and the related linear systems 
with fractional power of sparse positive definite (SPD) matrices
$\wcalAt^\alpha$ are discussed in Section~\ref{sec:3}.
Section~\ref{sec:4} is devoted to methods based on extensions of the underlying  
PDEs to domains in $\R^{d+1}$. Here we discuss two cases: 
extension to an elliptic problem in a semi-infinite cylinder, and
extension to a time-dependent problem.
Further, methods using integral representations of 
$\calA^{-\alpha}$ are considered in Section~\ref{sec:5}. The 
Balakrishnan integral and sinc-quadrature approximations are
surveyed first, followed by methods utilizing some alternative
integral formulas and quadratures.
The common in the methods from the last two sections is that they
solve numerically some reformulation of the fractional diffusion
problem. The BURA methods presented in Section~\ref{sec:6} follow 
a different approach. They approximate directly the inverse of the matrix
$\wcalAt^\alpha$. The best uniform rational approximation
of a properly defined scalar function on $[0,1]$ is used for this purpose. 
The computational efficiency is crucial in the case of
large-scale applications. This is the topic of the comparative 
analysis presented in Section~\ref{sec:7}, where issues related to
computational complexity and parallel scalability are discussed.  
Section~\ref{sec:8} is devoted to some challenges related to
numerical solution of time dependent space-fractional diffusion 
problems and coupled problems involving fractional diffusion 
operators.
Short concluding remarks are given at the end. \\

%\newpage
%%%%%%%%%%%%%%%%%%%%%%%%%%%%%%%%%%%%%%%%%%%%%%%%%%
\section{Spectral space-fractional diffusion problems}\label{sec:2}
%2 pages

\setcounter{section}{2}
\setcounter{equation}{0}\setcounter{theorem}{0}

\subsection{Spectral fractional powers of elliptic operators} 
\label{ss:problem}
Now we go to more general case of self-adjoint elliptic problems. Namely, 
we consider the following second order elliptic equation with homogeneous Dirichlet data:
\beq
\bal
 - \nabla \cdot( a(x) \nabla v(x)) + \q(x) v(x)&= f(x),\qquad \hbox{ for }
x\in \Omega,\\
v(x)&=0, \qquad \hbox{ for } x\in \partial \Omega.
\eal
\label{strong}
\eeq
Here $\Omega $ is a bounded domain in $\RR^d$,  $d\ge 1$.
We assume that $0<a_0 \le a(x) $, $a_0$ is a constant,  and $\q(x) \ge 0$ for $x\in \Omega$.
With the problem \eqref{strong} we associate  an elliptic operator defined in terms of the weak form 
of \eqref{strong}, namely, $v(x)$ is the unique function in  $V= H^1_0(\Omega)$ satisfying
\beq
a(v,\theta)= (f,\theta)\qquad \Forall \theta\in V.
\label{weak}
\eeq
Here
$$
a(w,\theta):=\int_\Omega  \Big (a(x) \nabla w(x) \cdot \nabla  \theta(x) + \q w(x) \theta(x) \Big )\, dx  
$$
and 
$$
(w,\theta):=\int_\Omega w(x) \theta(x)\, dx.
$$
For $f\in L^2(\Omega)$, \eqref{weak} defines a solution operator $\calT f :=v$.
Following \cite{kato}, we introduce an unbounded operator $\calA$ on $L^2(\Omega)$  as follows.
The operator $\calA$  with domain
 $$ D(\calA) = \{ \calT f\,:\, f\in L^2(\Omega)\}$$
 is defined by  
 \beq\label{def-A}
 \calA v=g \ \ \forall  v\in D(\calA), \ \ \mbox{ where } \ \  g\in L^2(\Omega) \ \ \mbox{ with } \ \  \calT g=v.  
 \eeq
 The operator $\calA$ is well defined as $\calT$ is injective.
 \begin{remark}\label{general}
 We note that the developed methods and algorithms are equally applicable to other than Dirichlet boundary conditions.
 For example one can assign Neumann or Robin boundary conditions or combination of all these. To avoid the 
 technical complications related to the case when the corresponding elliptic operator has its first eigenvalue zero,
 in all such cases we assume that the operator is positive definite, or equivalently, the corresponding bilinear form
 is coercive in the norm of the space $H^1(\Omega)$.
 \end{remark}
 \begin{remark}
 For those interested in the most general case of regularly accretive operators we refer to the paper
 of Bonito and Pasciak \cite{BP-17}.
 \end{remark}
The focus of this paper is numerical approximation and algorithm development for the
equation: 
\beq\label{frac-eq}
 \calA^\alpha u = f \quad  \mbox{with a solution} \quad u =\calA^{-\alpha} f.
\eeq
Here  $\calA^{-\alpha}=\calT^\alpha$ for $\alpha>0$
is defined by Dunford-Taylor integrals, \cite{kato}, which can be transformed  when
$\alpha\in (0,1)$, to the  Balakrishnan integral, e.g. \cite{balakrishnan}: for $f \in L^2(\Omega)$,
\beq\label{bal}
u=\calA^{-\alpha} f =\frac {\sin(\pi \alpha)} \pi 
\int_0^\infty \mu^{-\alpha} (\mu \calI +\calA)^{-1}f\, d\mu.
\eeq
This definition is sometimes referred to as the spectral definition of fractional powers.   
One can also use an equivalent definition through the 
expansion with respect to the eigenfunctions $\psi_j$ and the eigenvalues $\lambda_j$ of $\calA$, 
e.g. \cite{Acosta-Borth2017,LPGSGZMCMAK-18}:
\beq\label{eq:spectral-0}
\calA^{\alpha} u= \sum_{j=1}^\infty \lambda_j^{\alpha} (u,\psi_j) \psi_j \ \ \mbox{so that }
u = \sum_{j=1}^\infty \lambda_j^{-\alpha} (f,\psi_j) \psi_j. 
\eeq
Since the bilinear form $a(\cdot, \cdot)$ is symmetric on $ V \times V$ and $\calA$ is an unbounded 
operator we can show that $\lambda_j $ are real and positive and $\lim_{j \to \infty} \lambda_j =\infty$.

\begin{remark}
The operator $\calA$ defined by \eqref{def-A} preserves the positivity, that is,  
$\calA^{-1} f\ge 0$ whenever $f\ge 0$.
We note that by the maximum principle, if  $f \ge 0$  then $(\mu \calI+ \calA)^{-1} f \ge 0$ for $\mu\ge 0$,
and  then from \eqref{bal} we conclude that $\calA^{-\alpha}f \ge 0$.  In many applications, it is important that the 
corresponding approximations share this property.
\end{remark}

\begin{remark}
Another possible model of sub-diffusion reaction is given by the operator equation: find $u \in V$ s.t. 
\beq\label{eq:sub-reaction}
\calA^\alpha u + \qq u =f \quad \mbox{where} \quad \qq=const \ge 0.
\eeq
This problem could arise also by discretization of time-dependent sub-diffusion problems of the 
type $\frac{\partial u}{\partial t} + \calA^\alpha u= f(x,t)$ using a time stepping method with
some initial conditions, see Subsection \ref{sec:8.1}.
\end{remark}

\subsection{Regularity properties}\label{sec:2.2}
The regularity of the solution $u$ of the problem \eqref{frac-eq}  plays an essential role in devising and 
developing efficient numerical methods. The properties of the solution depend on the data $f$, the domain $\Omega$ and the
parameter $\alpha$. It is well known, that depending on these data the solution may develop singularities, boundary and/or 
internal layers that have to be captured by the numerical method.

The properties of the solution of \eqref{frac-eq} for the two basic definitions 
of fractional Laplacian  differ substantially. For example, the behavior near the boundary of the solution of the 
problem involving the {\it spectral fractional Laplacian}, e.g. \cite{CaffareliStinga16}, is
$$
\begin{aligned}
u(x) & \approx \text{dist}(x, \partial \Omega)^{2\alpha} +v(x)  \ \mbox{for} \ \alpha <1/2;\\
u(x) & \approx \text{dist}(x, \partial \Omega) + v(x) \ \ \mbox{for} \ \ \frac12 \le \alpha < 1,
\end{aligned}
$$
while the behavior of the corresponding problem involving 
{\it integral fractional Laplacian}, e.g. \cite{grubb2016,Ros_Oton_2014}, is
$$
u(x) \approx \text{dist}(x, \partial \Omega)^\alpha + v(x).
$$
Obviously, the low regularity of the solution near the boundary will lead to reduced order of convergence.

In this survey, we shall deal with the first definition, namely, spectral fractional Laplacian and its extension to more
general elliptic problem. For an extensive discussion of the numerical methods for the integral fractional Laplacian we 
refer to the review papers \cite{DU2020numerical, LPGSGZMCMAK-18, bonito2019numerical}.

%%%%%%%%%%%%%%%%%%%%%%%%%%%%%%%%%%%%%%%%%%%%%%%%%
\section{Discretization of the elliptic operator}\label{sec:3}

\subsection{Approximations of elliptic problems: main notations}
Here we shall give the main notations in discretizing the elliptic operator $ \calA $ by using
an $N$-dimensional finite element 
space $V_h$ or finite differences over a uniform mesh with $N$ points.  
In the case of finite element  approximation we get a 
symmetric and positive definite operator $\calAt: V_h \to V_h$, so that the approximation to \eqref{frac-eq} 
results in an operator equation $\calAt^{\alpha} u_h = f_h$ for a given $f_h \in V_h$ and the unknown $u_h \in V_h$. 
In the case of finite difference approximation 
we get a symmetric and positive
definite matrix $\wcalAt \in \RR^{N\times N}$ and a vector $\tilfh \in \RR^{N}$, so that the approximate solution  
$\tiluh \in \RR^{N}$ 
satisfies $\wcalAt^{\alpha} \tiluh = \tilfh$. 
These equations generate the so-called semi-discrete problems
\beq\label{semi}
\calAt^{\alpha} u_h =f_h \ \
( u_h = \calAt^{-\alpha} f_h ) \ \mbox{or/and} \ \  \wcalAt^{\alpha} \tiluh =  \tilfh \ \ ( \tiluh = \wcalAt^{-\alpha} \tilfh),
\eeq
where the fractional power is defined though the 
Balakrishnan integral formula \eqref{bal} or by \eqref{eq:spectral-0} with finite summation.
Below we give some particular examples of discretization using finite elements and finite differences.

\setcounter{section}{3}
\setcounter{equation}{0}\setcounter{theorem}{0}

\subsection{Finite element discretization}\label{sec:3.1}
The approximation in the finite element case is defined in terms of a
conforming finite dimensional space $V_h\subset V$ of piece-wise linear 
functions over a quasi-uniform partition of $\Omega$ 
into triangles or tetrahedrons.
Note that the  construction \eqref{bal}  of negative fractional powers
carries over to the finite dimensional case,  replacing $V$ and $L^2(\Omega)$ by
$V_h$ with $a(\cdot,\cdot)$ and $(\cdot,\cdot)$ unchanged.  

The discrete operator $\calAt$ is defined to be the inverse of
$\calTt:V_h\rightarrow V_h$ with
${\mathcal T}_h g_h:=v_h$  where $v_h\in V_h$ is the unique solution to
\beq
a(v_h,\theta_h) = (g_h,\theta_h),\Forall \theta_h\in V_h.
\label{T_h}
\eeq
The finite element approximation $ u_h \in V_h$ of $u$ is then given by
\beq
\calAt^{\alpha} u_h = \pi_h f, \ \ \mbox{or equivalently} \ \ \ 
u_h = \calAt^{-\alpha} \pi_h f:=\calAt^{-\alpha}f_h,
\label{fea}
\eeq
where $\pi_h$ denotes the $L^2(\Omega) $ projection into $V_h$.  In this
case, the dimension $N$ of the space $V_h$ equals the number
of (interior) degrees of freedom.  The operator $\calAt$ in the finite element
case is a map of $V_h$ into $V_h$ so that 
$\calAt v_h:=g_h$,  where $g_h\in V_h$ is the unique solution to \eqref{T_h}.

Let $\{\phi_j\}$ denote the standard ``nodal" basis of $V_h$.
In terms of this basis  $\calAt$ corresponds to the matrix
\beq\label{FEM-matrices}
\wcalAt = \wcalMt^{-1} \wcalSt, \ \  
\mbox{where} \ \ \wcalSt_{i,j} =a(\phi_i, \phi_j), \ \ \ \wcalMt_{i,j} =(\phi_i, \phi_j). 
\eeq
In the terminology of the finite element method, $\wcalMt$ and $ \wcalSt$ are the mass 
(consistent mass) and  stiffness matrices, respectively.

Obviously, if $\theta =\calAt \eta$ and ${\boldsymbol \theta},{\boldsymbol \eta} \in \R^N$ are the coefficient vectors
corresponding to $\theta,\eta\in V_h$, then $\boldsymbol \theta = \wcalAt
\boldsymbol \eta$.  Now, for the coefficient vector $\tilfh$ corresponding to $f_h=\pi_hf$ we have 
 $\tilfh = \wcalMt^{-1} \tilF$, where $\tilF$ is the vector with entries
$$
{\tilF}_j=(f,\phi_j), \qquad \hbox{ for }j=1,2,\ldots,N.
$$
Then using vector notation so that $\tiluh$ is the coefficient vector representing the solution $u_h$ 
through the nodal basis, we can write the finite element approximation of \eqref{strong} 
in the form of an algebraic system 
\beq \label{classic-FEM}
\wcalAt \tiluh = \wcalMt^{-1} \tilF \ \ \mbox{which implies} \ \ \wcalSt \tiluh = \tilF.
\eeq
Note that the matrix $\wcalSt$ is sparse while in general $\wcalAt$ is not.
However, when solving the standard diffusion problem \eqref{strong} one uses the sparse system $ \wcalSt \tiluh = \tilF$.

Consequently, the finite element approximation of the sub-diffusion problem \eqref{fea} becomes
\beq\label{mat-FEM}
\wcalMt \wcalAt^\alpha \tiluh = \tilF  \quad \hbox{or}\quad \tiluh=
\wcalAt^{-\alpha} \wcalMt^{-1} \tilF.
\eeq

We shall also introduce the finite element method with ``mass lumping" for two reasons.
First, it leads to positivity preserving fully discrete methods.
Second, it is well known that on uniform meshes lumped mass schemes for linear elements are equivalent
to the simplest finite difference approximations. This could be 
used to study the convergence of the finite difference method 
for solving the problem \eqref{frac-eq}, e.g., see \cite{harizanov2020JCP}.  

We introduce the lumped mass (discrete) inner product $ (\cdot,\cdot)_h$ on $V_h$ in the
following way (see, e.g. \cite[pp.~239--242]{Thomee2006})  for $d$-simplexes in $\RR^d$:
\beq\label{mass-lumping}
(z,v)_h = \frac{1}{d+1} \sum_{\tau \in {\mathcal T}_h } \sum_{i=1}^{d+1} |\tau|  z(P_i) v(P_i)
\ \ \mbox{and  } \ \  {\wcalMt}_h =\{ (\phi_i,\phi_k)_h\}_{i,k}^N.
\eeq
Here  $P_1, \dots, P_{d+1}$ are the vertexes of the $d$-simplex $\tau$ and $|\tau|$ is its 
$d$-dimensional measure.  The matrix $\wcalMt_h$ is called lumped mass matrix.
Simply, the ``lumped mass" inner product 
is defined by replacing the integrals determining the finite element mass matrix by local 
quadrature approximation, specifically, the quadrature defined by
summing values at the vertices of the $d$-simplex $\tau$ weighted by its measure $|\tau |$.   

In this case, we define $\calAt$ by 
$\calAt v_h:=g_h$  where $g_h\in V_h$ is the unique solution to
so that  $\calAt$ corresponds to the matrix 
\beq\label{A-lumped} 
\wcalAt = {\wcalMt}_h^{-1} \wcalSt. 
\eeq
Since ${\wcalMt}_h$ is diagonal matrix with
positive entries, the matrix $\wcalAt$ is sparse.
We also replace $\pi_h$ by $\calI_h$ so that the lumped mass semi-discrete
approximation is given by
\beq u_h = \calAt^{-\alpha} \calI_h f := \calAt^{-\alpha} f_h\quad \hbox{or} \quad   \bfu = 
\wcalAt^{-\alpha} \tilF.
\label{lumped-semi}
\eeq
Here $\tilF$ is the coefficient vector in the representation of the finite element interpolant $\calI_hf$ of $f$ 
with respect to the nodal basis in $V_h$.

\subsection{Finite difference discretization}\label{sec:3.2} 
In this case the approximation $\tiluh \in \R^N$  of $u$ is given by
\beq
\wcalAt^{\alpha} \tiluh = \tilfh, \ \ \mbox{or equivalently} \ \ \ 
\tiluh = %\calAt^{-\alpha} \cI_h f:=
\wcalAt^{-\alpha}  \tilfh, 
\label{fda}
\eeq
where $\wcalAt$ is an $N\times N$
matrix coming from  a finite difference approximation to the differential
operator appearing in \eqref{strong}, $\tiluh$ is the vector in $\R^N$ of the 
approximate solution at the interior $N$ grid points,  and 
$\tilfh \in \R^N$ denotes the vector of the values of  $f$ at the interior grid points.  
On a uniform mesh the matrix $\wcalAt$ is symmetric and positive definite. 

\paragraph{\it Example 1} We first consider the one-dimensional equation \eqref{strong} with
variable coefficient, namely, we study the following boundary value problem
$ - ( a(x) u^{\prime})^{\prime} =f(x),$ $ u(0)=0, \ u(1)=0, \ $ for $ 0<x<1$,
where $a(x)$ is uniformly positive function on $[0,1]$. On a uniform mesh $x_i =ih$, $i=0, \dots, N+1$,
$ h=1/(N+1)$,  we consider the three-point finite difference approximation of the second derivative
\begin{equation*}
\begin{split}
 (a(x_i) u'(x_i))' &
\approx %\frac{1}{\widetilde h_i} \left ( 
 \frac{1}{h}\left (a_{i+\frac12}\frac{u(x_{i+1}) - u(x_i)}{h} - a_{i-\frac12}\frac{u(x_i) - u(x_{i-1})}{h} \right ).
\end{split}
\end{equation*}
Here $a_{i-\frac12}=a(x_i - h/2)$ or  $a_{i - \frac12}=\frac{1}{h}\int_{x_{i-1}}^{x_i} a(x) dx $. Note that the former is 
the standard finite difference approximation obtained from 
the {\it balance method} (see, e.g. \cite[pp. 155--157]{samarskii2001theory}).

Then the finite difference approximation of the differential equation $- (a(x) u'(x))' =f(x)$ is given by the 
matrix equation \eqref{fda} with 
\begin{equation}\label{FD-matrix-1D-k} 
\wcalAt= \frac{1}{h^2} \left |
\begin{array}{ccccc} a_{\frac12}+ a_{\frac32}& -  a_{\frac32}  &&&\\
 -  a_{\frac32} &  a_{\frac32} +  a_{\frac52} & -  a_{\frac52}&&\\
\cdots &\cdots &\cdots &\cdots &\cdots\\ 
&&& - a_{N-\frac12} &  a_{N-\frac12} + a_{N +\frac12}\end{array}\right |. \ \ 
\end{equation}
The eigenvalues $\lambda_i$ of the matrix $\wcalAt $ are all real and positive and satisfy 
$$ 
4 \pi^2 \min_x a(x) \le \lambda_i \le  \frac{4}{h^2} \max_x a(x),  \ \  i=1, \dots, N.
$$

\paragraph{\it Example 2.}
The next example is for problem \eqref{strong}
on $\Omega=(0,1) \times (0,1)$ on a $(n+1) \times (n+1)$ square mesh with mesh-size $h=1/(n+1)$.
The standard 5-point stencil finite difference approximation of the Laplace operator gives the 
matrix  $\wcalAt \in \R^{N \times N}$, $N=n^2$, that has the
following block stricture (here $ \wcalAt_{i,i} \in \R^{n \times n}$, $i=1, \cdots, n $ 
and $\wcalIt_n$ is the identity matrix in $\R^n$) 
\begin{equation*}
\wcalAt=%\frac{1}{h^2}
\frac{1}{h^2} \left |\begin{array}{ccccc} \wcalAt_{1,1} & -\wcalIt_n &&&\\ -\wcalIt_n & \wcalAt_{2,2} & -\wcalIt_n 
&&\\\cdots &\cdots &\cdots &\cdots &\cdots\\ & -\wcalIt_n & \wcalAt_{i,i} & -\wcalIt_n &\\\cdots &\cdots &\cdots &\cdots 
&\cdots\\&&& -\wcalIt_n & \wcalAt_{n,n}\end{array}\right |,
\quad \wcalAt_{i,i}=\left |\begin{array}{cccc} 
4&-1&&\\-1&4&-1&\\\cdots&\cdots&\cdots&\cdots\\&-1&4&-1\\&&-1&4\end{array}\right |.
\end{equation*}

This matrix could be obtained by the finite element method applied to triangular meshes that 
generated on triangulations obtained by splitting each rectangle into two 
triangles (by connecting the lower left vertex with the upper right one)  
and using the ``lumped" mass inner product \eqref{mass-lumping}.
On a square mesh all diagonal elements of $ {\wcalMt}^{-1}_h$ are equal to 
$ h^{-2}=(n+1)^{2}$.
Then the finite element operator
$ {\calAt}: V_h \to V_h$ is defined as 
$ ({\calAt} u_h, v)_h = a(u_h,v) $ that results in the following matrix representation
$  \wcalAt= {\wcalMt}_h^{-1} {\wcalSt}$,
see \cite[Chapter 4, p.~203--205]{ciarlet2002}, where $ \wcalMt_h$ is the "lumped mass" matrix.

\begin{remark}\label{r:FD}
We note that on a uniform mesh with step-size $h=1/(n+1)$  the matrix  $\wcalAt$  
has the following extreme eigenvalues:
$$
\lambda_{min}=\lambda_{1,h} = 8(n+1)^2 \sin^2 \frac{\pi}{2(n+1)} \approx 2 \pi^2, $$ 
and  
$$
\lambda_{max} = \lambda_{n^2,h} = 8(n+1)^2 \sin^2 \frac{\pi n}{2(n+1)} \approx 8 (n+1)^2 = 8h^{-2}.
$$
\end{remark}
\paragraph{\it Example 3} 
We finally consider the one-dimensional equation \eqref{strong}, namely,
$-\Delta u := - u^{\prime\prime} =f(x),$ $ u(0)=0, \ u(1)=0, \ $ for $ 0<x<1$. On an arbitrary 
nonuniform grid $0=x_0 < x_1 < \dots < x_N < x_{N+1}=1$ we consider the three-point 
approximation of the second derivative
\begin{equation}\label{FD-matrix-1D} 
 %\calSt= \left  |
 \wcalSt = \left |
\begin{array}{ccccc} \frac{1}{h_1}+\frac{1}{h_2} & -\frac{1}{h_2} &&&\\
 - \frac{1}{h_2} & \frac{1}{h_2} + \frac{1}{h_3} & - \frac{1}{h_3}&&\\
\vdots &\vdots &\vdots &\vdots &\vdots\\ 
& -\frac{1}{h_i} &  \frac{1}{h_{i}} + \frac{1}{h_{i+1}} & -\frac{1}{h_{i+1}} &\\
\vdots &\vdots &\vdots &\vdots  &\vdots\\
&&& -\frac{1}{h_N} & \frac{1}{h_N} +\frac{1}{h_{N+1}} 
\end{array}\right |, \ \ 
\end{equation}
and $ \tilfh=[ \widetilde h_1 f(x_1), \widetilde h_2 f(x_2), \dots,  \widetilde h_{N+1} f(x_N) ]^T$,
  where $h_i=x_i -x_{i-1}$ and $\widetilde h_i = \frac12(h_{i+1}+h_i)$.

Note that this scheme is produced by the lumped mass finite element method with continuous 
piece-wise linear functions. Then $\wcalSt$ is  the stiffness matrix from linear finite elements
and $ \wcalMt_h =diag(\widetilde h_1, \widetilde h_1, \dots,\widetilde h_N )$
the lumped mass matrix so that  $\wcalAt  = \wcalMt_h^{-1} \wcalSt $.  
Then the finite difference  problem is $\calAt \tiluh =\tilfh$ and in this case $\wcalAt$ 
is symmetric in the inner product $\tiluh^T \wcalMt_h \tilvh : =(u_h, v_h)_h$.
On a nonuniform
mesh  in 1-dimension, this is well known, see, e.g. \cite[page 479]{samarskii2001theory}.

\subsection{Linear systems with fractional power of SPD matrices}\label{sec:3.3} 
The eigenvalues $\lambda_{j,h}$, $j=1, \dots, N$ of $ \calAt $ are real and positive and the eigenvectors 
$\psi_{j,h} \in V_h$, $j=1, \dots, N$ form a basis (ortho-normalized in 
$L^2$-inner product by $(\psi_{i,h},\psi_{j,h})=\delta_{ij}$,
where $\delta_{ij}$  equals $0$ if $j \not = i$ and is equal to $1$ for $j=i$).
The relation $\calAt \psi_{j,h}= \lambda_{j,h} \psi_{j,h}$ is equivalent to $\wcalSt \bfpsi_{j} = \lambda_{j,h} \wcalMt \bfpsi_j$,
where 
the entries of the vector columns  $\bfpsi_{j} \in \RR^N$ are the coefficients of $\psi_{j,h}$ represented through the 
nodal basis in $V_h$. Thus, orthonormality of the eigenvectors means $\bfpsi_{j}^T \wcalMt \bfpsi_{i} = \delta_{ij}$.

Now after introducing the matrix $\Lambda=diag(\lambda_{1,h}, \dots, \lambda_{N,h})$ and the 
matrix $\Psi \in \RR^{N\times N}$ with columns the eigenvectors $ \bfpsi_{j}$ we have 
the following  representation of the solution \eqref{fda}: $\tiluh = \Psi^T  \Lambda^{-\alpha} \Psi \ \tilfh$.
This could lead to quite 
fast and reliable solution method, especially, if FFT is applicable. Unfortunately,  the application of such method 
will be limited to 
Laplace operator and rectangular domains and uniform meshes.
For general domains and variable coefficients, computing the whole spectrum  from 
 $\wcalSt \bfpsi_{j} = \lambda_{j,h} \wcalMt \bfpsi_j$, $j=1, \dots, N$, 
 would be a prohibitively expensive procedure for large $N$. 

Thus, we need a method for approximately solving 
\eqref{fea}. Such methods will be called {\it  fully discrete} methods reflecting the situation that we first
define the discrete fractional order to the elliptic problem that leads to a system of linear equations and 
then we apply an approximate method for solving that problem. Here we survey a number of such methods.

As shown recently in \cite{HOFREITHER2019}, 
these methods, though entirely different, are interrelated and all seem to involve certain rational approximation of
the fractional powers of the underlying elliptic operator. As such, from mathematical point of view,
those based on the best uniform rational approximation should be the best. However, one should realize that BURA
methods involve application of the Remez method for finding the best uniform rational approximation,
\cite{PGMASA1987,Driscoll2014}, a numerical algorithm for solving certain min-max problem that is highly 
nonlinear and sensitive to the precision of the computer arithmetic. For example, in \cite{varga1992some}
 the best uniform rational approximation of $z^\alpha$ for six values of $\alpha  \in (0,1)$ 
are reported for degree $k\le 30$ by using computer arithmetic with $200$ significant digits.

%%%%%%%%%%%%%%%%%%%%%%%%%%%%%%%%%%%%%%%%%%%%%%%%%
\section{Methods based on extensions to PDEs  in domains in $\RR^{d+1}$}\label{sec:4}

\setcounter{section}{4}
\setcounter{equation}{0}\setcounter{theorem}{0}

\subsection{Extension to elliptic problem in a semi-infinite cylinder}
\label{sec:4.1}
We shall demonstrate this approach in the case of fractional Laplacian and the concept of 
“Neumann to Dirichlet” map proposed in  \cite{caffarelli2007extend} to study the existence, uniqueness
and  the regularity of the solution of \eqref{fea}.  
Namely, the solution of fractional Laplacian problem is obtained from the relation  
$u(x) = U(x, 0)$ 
where $U : \Omega\times\RR_+ \rightarrow \RR$ is a solution of the equation
$$
-div\left ( y^{1-2\alpha}\nabla U(x,y)\right ) = 0, ~~~ 
(x,y)\in \mathcal C= \{\Omega\times \RR_+ \}.
$$
Here $U(\cdot,y)$ satisfies the boundary conditions in \eqref{strong} and 
in addition 
$$
\lim_{y\rightarrow\infty} U(x,y) = 0
\ \ 
\mbox{as well as}  \ \ 
\lim_{y\rightarrow 0^+} \left (- y^{1-2\alpha} U_y(x,y)\right ) = f(x), 
~~~ x\in\Omega.
$$
The variational formulation of this equation is  posed in some weighted  Sobolev space, \cite{caffarelli2007extend},
$$
\dot H^1_L({\mathcal C},y^\alpha)=\{ w \in H^1({\mathcal C}, y^\alpha):  w=0 \ \mbox{on } \ \partial_L {\mathcal C} \},
$$ 
where $ \partial_L {\mathcal C} $ is the lateral surface of the infinite domain ${\mathcal C}$ and
$$ 
H^1({\mathcal C},y^\alpha)=\{ w: \int_{\mathcal C}y^{1-2\alpha} ( |\nabla w|^2 + w^2) dxdy < \infty\}. 
$$
Then, one seeks a solution $U \in \dot H^1_L({\mathcal C},y^\alpha)$ that satisfies the 
following  integral identity (to make it simpler we assume that $f \in L^2(\Omega)$):
$$
\int_{\mathcal C} y^{1-2 \alpha} \nabla U \cdot \nabla \phi dx dy = d_\alpha \int_\Omega f \phi  dx, \ \ \  
\forall \phi \in \dot H^1_L({\mathcal C},y^\alpha).
$$
Here $\nabla U$ is the gradient of $U$ in the variables $(x,y)$ and $d_\alpha$ is a normalizing constant 
$d_\alpha=2^{1-2\alpha}\Gamma(1-\alpha)/\Gamma(\alpha)$, see, e.g. \cite[formula (2.26)]{nochetto2016pde}.

The finite element 
approximation, proposed and studied in \cite{nochetto2016pde,CNOSalgado-2016}, 
uses the rapid decay of the solution $v(x,y)$ in the 
$y$ direction, thus enabling truncation of the semi-infinite
cylinder to a bounded domain of modest size, namely ${\mathcal C}_Y = \Omega \times (0, Y)$. 
Then the finite element approximation $U_h(x,y)$ of $U(x,y)$ is sought as a solution of the weak form in the finite dimensional subspace
$X_h= V_h \times W_h$, where $W_h$ is the set of piece-wise linear functions on a partition of $(0,Y)$ 
and $W_h \subset \{ w \in H^1(0,Y), w(Y)=0\}$. If the dimension of $W_h$ is $M$, then the count of all mesh points is $NM$.
An almost optimal with respect to the number of the degrees of freedom (with a log-factor) rate of decay of the error 
$
\| u (x)- U_h(x,0)\|_{H^\alpha(\Omega)}
$
has been established in   
\cite[see, Theorem 5.4 and Remark 5.5]{nochetto2016pde}. The authors use delicate and sharp analysis, proper  choice of $Y$, 
and graded near $Y=0$ meshes.
Further in \cite{CNOSalgado-2016} an efficient multilevel method based on the Xu-Zikatanov identity \cite{XuZ2002}
has been proposed, studied and tested. 

Recently, Hofreither, \cite{HOFREITHER2019}, made an interesting interpretation of this method by rewriting it 
in the following way. First, using separation of variables, one introduces the following eigenvalue problem in direction $y$ 
(see, \cite{HOFREITHER2019} and also \cite[problem (2.25)]{nochetto2016pde} with different normalization):
\begin{equation}\label{EVP-W}
-(y^{1-2\alpha}\psi_k^{\prime})^ {\prime}  = \mu_k y^{1-2\alpha} \psi_k, \ \     0< y < \infty , \ 
\psi_k(0)=1, \ \lim_{y \to \infty} \psi_k(y)=0.
\end{equation}

The approximation of this problem in the truncated  interval $(0,Y)$ and boundary condition $\psi_k(Y)=0$
on the finite element space $ W_h$  of dimension $M$ produces   
the eigenpairs $(\mu_{k,h}, \psi_{k,h}(y))$, $k=1, \dots, M$.  Then by separation of variables
we get the following representation of the solution 
$$
U_h(x,y)=d_\alpha\sum_{k=1}^M \sum_{j=1}^N \psi_{k,h}(y) \psi_{j,h}(x) \frac{\psi_{k,h}(0)}{\mu_{k,h} + \lambda_{j,h}} (f_h,  \psi_{j,h}).
$$
Therefore, 
\begin{equation}\label{H-RA}
U_h(x,0) = \Psi^T  r(\Lambda) \Psi f_h, \ \mbox{where} \  r(z)=d_\alpha\sum_{k=1}^M \frac{\psi_{k,h}(0)^2}{\mu_{k,h} + z}
\end{equation}
with $\Lambda=diag(\lambda_{1,h}, \dots, \lambda_{N,h})$. Since $u_h = \calA^{-\alpha} f_h = \Psi^T \Lambda^{-\alpha} \Psi$ 
and $U_h(x,0)$ approximates $u_h $, obviously $r(z)$ is a rational approximation of $z^{-\alpha}$.

The approximation (\ref{H-RA}) can be expressed in terms of tensor products.   
The analysis in \cite{HOFREITHER2019} fully decouples the error 
in the extended direction $y$ from the error in the spatial variable $x$.
This framework result may be used for further elaboration of estimates 
in $L^2(\Omega)$, as opposed to the error estimates in \cite{nochetto2016pde}, 
which are in weighted fractional Sobolev spaces. 

The results in \cite{GHH-19,HOFREITHER2019} show certain advantages of using 
discretization of higher order to define the space $W_h$.  Some numerically 
computed convergence rates of system solves versus the dimension
$M$ of $W_h$ are given in Table \ref{E&LvC}. 
\renewcommand{\arraystretch}{1.2}
\begin{table}[ht]
\caption{ \footnotesize Convergence rate of system solves: discretization of $W$ 
with linear splines from $C^0[0,Y]$ versus cubic splines from $C^2[0,Y]$}
\label{E&LvC}
\centering
\begin{tabular}{|c|c|c|}
\hline
$\alpha$ & linear FEM & cubic splines\\
\hline
0.25 & $M^{-2}$   & $M^{-6}$ \\
0.50 & $M^{-2}$   & $M^{-6}$ \\
0.75 & $M^{-2.5}$ & $M^{-3}$ \\
\hline
\end{tabular}
\end{table} 
In opposite to the cases $\alpha\in\{0,25,0,50\}$, the rates 
for $\alpha=0.75$ are rather closer, which needs some more involved analysis. 

Now, we briefly comment the numerical stability of the eigenvalue 
problem (\ref{EVP-W}), where 
$M$ is supposed to be not very large. However, depending on the space
$W_h$ and even stronger depending on the value of $\alpha$, it may become
very ill conditioned. The loss of accuracy of the numerically computed
spectrum $(\mu_k,\psi_k)$ may practically destroy the accuracy of the
rational approximation (\ref{H-RA}). To stabilize the computations, a 
simple regularization procedure is proposed in \cite{HOFREITHER2019}.
In any case, one has to be careful at this point.
Some related issues are discussed in \cite{SXK-17}, where spectral FEM is 
applied to the fractional diffusion problem. As noted there, a loss of accuracy 
is caused if the eigenfunctions are not perfectly orthogonal. To deal with 
this, a weighted Gram-Schmidt orthogonalization is applied 
resulting to a significant improvement in the spectral FEM accuracy.

\subsection{Extension to time-dependent problem}\label{sec:4.2}

This is another method based on seeking a function $\Uhxt(x,t)$ on
the extended domain $\Omega \times (0,1)$. 
This approach is based on the following observation of 
Vabishchevich, \cite{Vabishchevich14,Vabishchevich15}: if 
$\calBt=\calAt -\delta\calIt$ with $\delta \in (0,  \lambda_{1,h}]$,
and $\Uhxt(t) %\tilde{u}_{h}(t)
\in V_h$ is the solution of the initial value problem
\begin{equation}\label{extend-eq-uh}
(\delta \calIt+t \calBt )\partial_t \Uhxt (t) + \alpha \calBt \Uhxt(t) =0 \ \  t \in (0,1]\  \
\Uhxt(0)=\delta^{-\alpha}  f_h,
\end{equation}
then $u_h=\Uhxt(1)$. 
Thus, the solution of the original problem \eqref{fea} is sought in the cylinder $(x,t) \in \{ \Omega \times (0,1) \}$
as a solution to a homogeneous pseudo-parabolic equation with initial data $\delta^{-\alpha}f_h$. The value of the solution 
at the final time $t=1$ represents the solution of \eqref{fea}.
In  \cite{Vabishchevich14,Vabishchevich15} Vabishchevich 
proposed and studied various two-level  schemes and showed optimal convergence rates for 
sufficiently smooth solutions $ \Uhxt(t)$  (with respect to time $t$). 
As outlined in \cite{DLP-Pade}, there is substantial difference in the smoothing properties of the 
pseudo-parabolic operator of \eqref{extend-eq-uh} and the one associated with the standard parabolic equation.
This leads to a completely different regularity pick-up of the solution of equation \eqref{extend-eq-uh} from the 
data $f$ compared with  the standard parabolic problem, which shows exponential decay of the solution,
see below Remark \ref{exponent}.  In short, to have smooth solution $ \Uhxt(t)$
one needs to assume high regularity and/or compatibility conditions of the right-hand-side $f$ with the 
boundary conditions. Problems with data that do not satisfy any of these are called problems with non-smooth data.
Following the original idea from \cite{Vabishchevich15}, 
various time-stepping schemes for solving the parabolic problem \eqref{extend-eq-uh} 
have been developed and studied. They all propose improvements of the original algorithm 
in making it more efficient 
and/or more general, e.g.  \cite{CV-19,CV-20,DLP-Pade,Vabishchevich18}.

Here we shall present a method based on this approach to problems with  non-smooth data,  e.g. \cite{DLP-Pade}
The discretization scheme for \eqref{extend-eq-uh} uses  
geometrically refined  near the origin mesh and Pad\'e approximations
of the function $(1+x)^{-\alpha}$ with rational functions of the type $P_m(x)/Q_m(x)$. 
A construction of the time-stepping mesh with a rigorous analysis when $m \to \infty$ is presented in \cite{DLP-Pade}
which is explained briefly below.

Using the spectrum 
of $\calAt$  we express $\Uhxt(t) \in V_h$ and $f_h \in V_h$ through the  expansions
\begin{equation*}
\Uhxt(t)=\sum_{j=1}^{\D} \hat{U}_{j,h}(t) \psi_{j,h}\qquad \mbox{and}\qquad {f}_{h}=\sum_{j=1}^{\D} \hat{f}_{j,h} \psi_{j,h}.
\end{equation*}
Substituting these  into  \eqref{extend-eq-uh} we get that $\hat{U}_{j,h}(t)$, $0 < t  \le 1$,  satisfies 
\begin{equation}\label{u-hj}
\hat{U}_{j,h}(t) = \frac{\hat{f}_{j,h}}{(\delta + t(\lambda_{j,h}-\delta))^{\alpha}}
\quad  \mbox{which implies } \ \ 
\hat{U}_{j,h}(1)=\lambda_{j,h}^{-\alpha} \hat{f}_{j,h}.
\end{equation}
Thus, we get
$
\Uhxt(1)=\sum_{j=1}^{\D} \lambda_{j,h}^{-\alpha} \hat{f}_{j,h}\psi_{j,h}=\calAt^{-\alpha}f_h = u_h.
$

\begin{remark}\label{exponent}
Note that the solution of the standard parabolic problem $\partial_t U_h + \calAt U_h=0  $ has an expansion
with respect to the eigenfunctions of $\calAt$ so that $ {U}_{h,j}(1)=\sum_{j=1}^N e^{-\lambda_{j,h}} \hat{f}_{j,h} \psi_{j,h}$,
which shows an exponential decay (with respect to the eigenvalues) of the initial data.
\end{remark}
 
Now we present a generalization and an improvement of  Vabishchevich method proposed in  \cite{DLP-Pade}.
The improvement is due to the use of a  diagonal Pad\'e approximation of $(1+z)^{-\alpha}$ 
for  $z\in [0,1]$:
\begin{equation}\label{pade-0}
(1+z)^{-\alpha}=\frac{P_m(z)}{Q_m(z)}+\epsilon_m(z):=r_m(z)+\epsilon_m(z),
\end{equation}
where $m \in \mathbb{N}^+$ and $P_m$, $Q_m$ are polynomials of degree $m$. 

Then for a given temporal mesh $0=t_0<t_1<\cdots<t_{M+1}=1$ we
introduce a discretization scheme through the following recursion
for the approximation $U_{l,h}$ of $  \Uhxt(t_l)$
\begin{equation}\label{rec-extend-Uh}
 {U}_{h, l+1} =r_m\left(k_l \calBt (\delta \calIt+t_l\calBt)^{-1} \right) {U}_{h,l},\quad l=0,1,\cdots, M,
\end{equation}
with step-size  
$k_l=t_{l+1}-t_l$, $l=0,1,\cdots, M$ 
and $ {U}_{h,0}=\Uhxt(0)=\delta^{-\alpha}f_h$. We will take ${U}_{h,M+1}$ as an  
approximation of $\Uhxt(1)$.  Note, that $m=1$ will produce  Crank-Nicolson scheme,
advocated in  \cite{CV-19,Vabishchevich18}.

The efficiency of the method will depend substantially on the choice of the mesh-points $t_n$.
This was discussed in details in  
\cite{DLP-Pade}, where the following  two-level construction of such meshes has been proposed:
(1) first introduce a geometrically  refined (near zero) mesh using the points $\ttt_0=0$, $\ttt_l = 2^{l-1-L}$, $l=1, \dots, L+1$,
and  (2) divide each subinterval $(\ttt_l, \ttt_{l+1})$ into $n$ equal subintervals. Now, if
$\tkk_l=(\ttt_{l+1}-\ttt_l)/n$, we get the following set of total $M+1$ (with $M=n(L+1)$) points in time direction
$$
0< \tkk_1 < \cdots < n \tkk_1=\ttt_1 < \ttt_1 + \tkk_2<   \ttt_1 + 2\tkk_2< \cdots  < \ttt_2< \cdots < \ttt_{L+1}=1.
$$
After renaming this mesh as 
 $0=t_0<t_1<\cdots<t_{M+1}=1$ we apply the approximation scheme \eqref{rec-extend-Uh}.
 This scheme was  studied and numerically tested in  \cite{DLP-Pade}  by considering that $m$ is fixed, say $m=1,2,3$, 
 while letting $n \to \infty$.
In \cite{DLP-Pade} the authors proved and experimentally confirmed algebraic convergence (while 
 refining the mesh in time and keeping $m$ fixed), namely,
 $$
 \|U_{h,M+1}  -  \calAt^{-\alpha}{f}_{h} \| \le C(\alpha, \delta)  n^{-2m} \| f_h \|,  \ \  n \ge 2.
 $$
As shown in \cite{DLP-Pade}, the method requires  $ n m (L+1)$ system solves of the type
$ (\delta \calIt+t_l\calBt) v = w$.

%%%%%%%%%%%%%%%%%%%%%%%%%%%%%%%%%%%%%%%%%%%%%%%%%
\section{Methods based on integral representation of $\calAt^{-\alpha}$}\label{sec:5}

\setcounter{section}{5}
\setcounter{equation}{0}\setcounter{theorem}{0}

\subsection{Approximation of Balakrishnan integral by sinc quatratures}
\label{sec:5.1}
This type of methods have been developed, theoretically justified and  practically tested in a series of papers by
Bonito and Pasciak,  \cite{bonito2019numerical,bonito2019sinc,BP-17}. With a simple change of the variable 
$\mu =e^y$ in \eqref{bal} and replacing $\calA$ by $\calAt$ we get the following semi-discrete solution
$$
u_h=\calAt^{-\alpha} f_h =\frac {\sin(\pi \alpha)} \pi \int_{-\infty}^\infty e^{(1-\alpha)y} (e^y \calIt +\calAt)^{-1}f_h\, dy.
$$
The proposed methods are based on a truncation of the integral and an application of 
a proper quadrature formula, \cite{BP-17, bonito2019sinc}. Here we shall present the most popular and the best (in 
terms of accuracy and smoothness requirements of the data $f$) method, given in \cite{bonito2019sinc}.

It is based upon a selection of a positive quadrature step $\kk$ and quadrature nodes $y_l=l\kk$ with two 
integers $M$ and $N$  that constitute the sinc approximation 
of the truncated Balakrishnan integral \eqref{bal} over $(-M\kk, N \kk)$  so that:
$$
u_{h,\kk}=Q_\kk^{-\alpha}(\calAt) f_h,
$$
where
\beq\label{BP-Q}
Q_\kk^{-\alpha}(\calAt) := \frac {\kk \sin(\pi \alpha)} \pi 
\sum_{l=-M}^N e^{(1-\alpha)y_l} (e^{y_l } \calIt +\calAt)^{-1}.
\eeq
The quadrature step $\kk$ is a real number and the integers $N$ and $M$ are taken to be 
of order $1/\kk^2$, cf. \cite[Remark 3.1]{BP-17,bonito2019sinc}. 
The error analysis is done by careful estimation of the quadrature error 
$$
\begin{aligned}
\calA^{-\alpha} f  - Q^{-\alpha}_\kk(\calA) f & = \int_{-\infty}^\infty F(y,f) dy - \kk \sum_{l=-\infty}^\infty F(l\kk) f  \\
      &  + \kk  \sum_{l <-M} F(l\kk) f  + \kk \sum_{l>N} F(l\kk) f, 
\end{aligned}
$$
with $F(y) =  e^{(1-\alpha)y} (e^{y} \calIt +\calAt)^{-1} $.
The last two terms  have different behavior with respect to $\alpha$ so having 
two different $M$ and 
$N$ allows to balance the three errors: the quadrature error and two errors due to truncating the infinite integral. 
As shown in \cite[Remark 7.3]{BP-17}  the choice
\begin{equation}\label{MN-choice}
M=\left \lceil \frac{\pi^2}{4 \alpha \kk^2} \right \rceil  \ \ \mbox{and}  \ \ N=\left \lceil \frac{\pi^2}{4 (1-\alpha) \kk^2} \right \rceil
\end{equation}
balances these errors.
Theoretically, this scheme has exponential rate of convergence as $\kk \to 0$. A simplified form of Theorem 4.2 in  \cite{bonito2019sinc}
gives the following estimate for $ u_h - u_{h,\kk}$:
$$
\| u_h - u_{h,\kk} \|_{L^2(\Omega)} \le C e^{-\pi^2/(2\kk)} \| f\|_{L^2(\Omega)}.
$$
In our numerical tests we call this scheme the $\kk$-Q-method. 
 We stress that due to the choice \eqref{MN-choice} of different quadrature points $M$ and $N$
 for the negative and positive semi-axis, this method is robust with respect to $\alpha \in (0,1)$, while some methods
 deteriorate substantially for $\alpha$ close to $0$. 
 
 A simplified version with $N=M \sim  1/\kk^2$, 
 though less efficient, is also used.  
The paper \cite{bonito2019sinc}  contains a number of other estimates, like convergence in $H^r(\Omega)$
for: (1) $\alpha > r/2$ and $f \in L^2(\Omega)$  and (2)  $\alpha \le r/2$  and $f \in H^{r -2\alpha+\epsilon}(\Omega)$,
$\epsilon >0$ (see, \cite[Assumption 4.1 and Theorem 4.2]{bonito2019sinc}. 
The choice \eqref{MN-choice} achieves an exponential rate with respect to the number of quadrature nodes $N$.
For each $y_l$ one solves the system  $  (e^{y_l } \calIt +\calAt) w_l=f_h$, $l=-M, \dots, N$,
which results in solving $N+M+1$ systems.

\subsection{An alternative method based on Gauss-Jacobi quadratures}
\label{sec:5.2}
After a change of variable $\xi=\tau \frac{1-\mu}{1+\mu}$, $\tau >0$ in \eqref{bal},  we get 
\beq\label{Ace-Nov}
\calAt^{-\alpha} =C_\alpha \tau^{1-\alpha}
\int_{-1}^1 (1- \xi)^{-\alpha} (1+\xi)^{\alpha -2} \left ( \tau \frac{1-\xi}{1+\xi} \calIt +\calAt \right )^{-1} \, d \xi,
\eeq
with $C_\alpha=  {2\sin(\pi \alpha) } /\pi$.
To approximate this integral a 
$k$-point Gauss-Jacobi rule with respect to the 
weight $\omega(\xi)=(1-\xi)^{-\alpha}(1+\xi)^{\alpha-1} $ has been proposed and studied in  \cite{AN-18,AN-19},
see also \cite{AN-20}:
\begin{equation}\label{Gauss-J}
 \calAt^{-\alpha}  \approx Q_k^{-\alpha} (\calAt):= \sum_{j=1}^k \gamma_j(\eta_j \calIt + \calAt)^{-1},
\end{equation}
with
$$
\gamma_j =\frac {2\sin(\pi \alpha) \tau^{1-\alpha}}{\pi} \frac{w_j}{1+\theta_j}, \ \ \ 
\eta_j= \frac{\kappa(1-\theta_j)}{1+ \theta_j}.
$$
Here $\omega_j$ and $\theta_j$, $j=1, \dots,k$, are, respectively,  the weights and the nodes of the Gauss-Jacobi quadrature.
The choice of $\tau$ is critical for the quality of the approximation 
of  \eqref{Ace-Nov} by \eqref{Gauss-J}. As shown in \cite{AN-17}, for large $k$ the optimal choice is
$ \tau =\sqrt{\lambda_{1,h} \lambda_{N,h}} $. The best practical choice of $\tau $ is provided in 
\cite[Proposition 4, formula (32)]{AN-17}.  The error analysis of the method, as presented in \cite{AN-17},
relies on the relation between the Gauss-Jacobi quadrature error
and the Pad\'e approximation 
$P_{k-1}(1-\zz)/Q_k(1-\zz)$ of $(1-\zz)^{-\alpha}$ on the interval $(-1,1)$. It is expressed thorough the Pad\'e
approximation error
$$
E_{k-1,k}(1-\zz) := (1-\zz)^{-\alpha} -   P_{k-1}(1-\zz)/Q_k(1-\zz),
$$
so that 
\begin{equation}\label{Aceto-error}
\|   \calAt^{-\alpha}  -  Q_k^{-\alpha} (\calAt)  \| \le \max_{\lambda \in [\lambda_{1,h}, \lambda_{N,h} ]}
\tau^{-\alpha} E_{k-1,k} \left (\lambda/\tau \right ).
\end{equation}
The optimal choice of $\tau$ is obtained by minimizing the right hand side of \eqref{Aceto-error} for $\tau \in (0, \infty)$. 
The minimization problem is solved approximately and the optimal parameter $\tau$ is shown to depend
on $\lambda_{1,h}$, $\lambda_{N,h}$, and $k$. The optimal choice of $\tau$ gives an  
asymptotic in $k$ error bound, \cite[Theorem 3]{AN-17}:
$$
\|   \calAt^{-\alpha}  -  Q_k^{-\alpha} (\calAt)  \| \le C \lambda_{N,h}^{-\alpha/2} e^{-ck/\sqrt[4]{\lambda_{N,h}}}.
$$
For a fixed mesh the error of this method shows exponential decay.
However, for fixed $k$ the factor $ e^{-ck/\sqrt[4]{\lambda_{N,h}}}$ tends to 1 when the 
mesh size goes to zero, thus the error deteriorates.  At the same time the first factor $ \lambda_{N,h}^{-\alpha/2} $ tends to zero, so the 
convergence is always ensured. The numerical experiments provided in \cite{AN-17} illustrate adequately the error behavior. 

In the same spirit, but using different idea, is the approach proposed by Vabishchevich in \cite{V-20}, 
based on change of the variable  $\mu = \xi(1- \xi)^\sigma$, $\sigma >0$,
so that
$$
\calAt^{-\alpha}= \frac{\sin(\alpha \pi)}{ \alpha \pi}\int_{0}^1 (1-\xi)^{\sigma\frac{1-\alpha}{\alpha}-1}  
(1+(\sigma-1)\xi) \left ( (1-\xi)^\frac{\sigma}{\alpha} \calIt + \xi^{\frac{1}{\alpha}} \calAt \right )^{-1} d\xi.
$$
If $\calAt$ is properly normalized so that $\lambda_{1,h}=1$ 
then
$$
F(\xi, \zz) = (1-\xi)^{\sigma\frac{1-\alpha}{\alpha}-1}  
(1+(\sigma-1)\xi) \left ( (1-\xi)^\frac{\sigma}{\alpha}  + \xi^{\frac{1}{\alpha}} \zz \right )^{-1}, \ 1 \le \zz < \infty,
$$
considered as a function of $\xi$ does not have singularities in $0 \le \xi \le 1$. Moreover, 
for  $\sigma$ sufficiently large the function has continuous derivatives of high order.
Other transformations are also possible, see, e.g.   \cite[e.g. Formulas (22), (23)]{V-20}.
Here we shall present the main idea, while the interested reader can find all details 
and relevant numerical experiments in \cite{V-20}.
If one chooses properly $\sigma$ 
then any standard 
composite quadrature rule with $M$ subintervals, e.g. composite trapezoidal or Simpson rules,  is applicable.
 For example, e.g. \cite{V-20}, the choice $\sigma =\kappa \max(\alpha^{-1}, (1-\alpha)^{-1})$
allows to use the  trapezoidal rule for $\kappa=2$ and the Simpson rule for $\kappa=4$. 
The theoretical estimates
of the error in terms of the data is to be developed yet, but the numerical experiments provided
in \cite{V-20} are quite promising. A possible downside of the method could  be an error bound 
that involves the norm of  $\calAt^\kappa f_h$. This will indicate a lost of the quality of 
the approximation for non-smooth data and in general, loss of robustness.

\subsection{Conclusions on integral formulas and quadratures}
\label{sec:5.3}
The discussed in this Section methods lead to algorithms that produce a particular rational approximation 
of $z^{-\alpha}$ and consequently produce an algorithm that requires solution of algebraic problems 
of the type $(\calAt + c\calIt)w=v$, $c>0$, which number is equal to the number of the quadrature points.
It is quite obvious that $ Q_\kk^{-\alpha}(\calAt) $  of Bonito-Pasciak method \eqref{BP-Q}  
and $Q_k^{-\alpha}(\calAt)$ of Aceto-Novaty method \eqref{Gauss-J} are rational functions of $\calAt$.
It is less obvious for the 
method \eqref{rec-extend-Uh} based on the pseudo-parabolic equation, \cite{CV-19,CV-20,DLP-Pade} that it also generates 
a rational function. Indeed, since  
$\calBt = \calAt- \delta \calIt$ we can rewrite Vabishchevich method \eqref{rec-extend-Uh} 
in the following form
$$
U_{h,M+1} = \prod_{l=0}^M  r_m\left(k_l (\calAt -\delta \calIt )(\delta (1-t_l) \calIt+t_l\calAt)^{-1} \right) {U}_{h,0}.
$$
Obviously, the operator $ \delta^{-\alpha}\prod_{l=0}^M  r_m\left(k_l (\calAt -\delta \calIt )(\delta (1-t_l) \calIt+t_l\calAt)^{-1} \right) $ 
advocated in  \cite{DLP-Pade,Vabishchevich14,Vabishchevich15} is a rational approximation of $\calAt^{-\alpha}$. 
Thus,  solving numerically the pseudo-parabolic equation \eqref{extend-eq-uh} based on various 
time-stepping strategies proposed, studied and tested in  \cite{CV-19,CV-20,Vabishchevich18},
could be interpreted as designing particular rational approximation of $\calAt^{-\alpha}$. 
 
%%%%%%%%%%%%%%%%%%%%%%%%%%%%%%%%%%%%%%%%%%%%%%%%%
\section{Methods based on the best uniform rational approximations}
\label{sec:6}

\setcounter{section}{6}
\setcounter{equation}{0}\setcounter{theorem}{0}

\subsection{Best uniform rational approximation of $\zz^\alpha$ on $[0,1]$  (BURA)}
\label{sec:6.1}

 In order to 
 use the known results for the approximation theory, we first rewrite the solution 
 of the \eqref{fda} in the form
 \beq\label{eq:sol-function}
 u_h = \lambda_{1,h}^{-\alpha} (\lambda_{1,h} \calAt^{-1})^\alpha f_h.
 \eeq
 The scaling by $\lambda_{1,h}$ maps the eigenvalues of $ \lambda_{1,h} \calAt^{-1}$ to 
 $(\lambda_{1,h}/\lambda_{N,h},1] :=(\delta, 1]  \subset (0,1]$. Here $\delta=\lambda_{1,h}/\lambda_{N,h}$ is a small 
 positive number.

Now we consider BURA along the diagonal of
the Walsh table and take $\rat k$ to be the set of rational functions 
$$
\rat k= \bigl\{ r(\zz): r(\zz)=P_k(\zz)/Q_k(\zz),   
\   P_k \in {\mathcal P}_k, \mbox{ and  } \ Q_k \in {\mathcal P}_k, \ \mbox{monic}  \bigr\}
$$
with  ${\mathcal P}_k$ set of algebraic polynomials of degree $k$. To find an approximation to 
\eqref{eq:sol-function} we introduce the  
best uniform rational approximation (BURA) $r_{\delta,\alpha,k}(\zz) $ of $z^\alpha$ on $[\delta,1]$
$$
r_{\delta,\alpha,k}(\zz) := \argmin_{s(\zz)\in \rat k}\,   \sup_{\zz \in [\delta,1]} |  s(\zz) - \zz^{\alpha} |.
$$
It is quite appealing to get rid of $\delta=\lambda_{1,h}/\lambda_{N,h}$ by 
using the best uniform rational approximation $ r_{\alpha,k}(\zz)$ of $\zz^\alpha$ on the whole interval $[0,1]$, namely
\begin{equation}\label{bura}
 r_{\alpha,k}(\zz) := \argmin_{s(\zz)\in \rat k}\,  
 \max_{\zz \in [0,1]} |  s(\zz) - \zz^{\alpha} | = \argmin_{s(\zz) \in \rat k}  \| s(\zz) - \zz^{\alpha} \|_{L^\infty(0,1)}.
\end{equation}
\begin{remark}\label{delta-BURA}
It is obvious that 
$$ \|r_{\delta, \alpha,k} (\zz) - \zz^{\alpha} \|_{L^\infty[\delta,1]} < \|r_{\alpha,k} (\zz) - \zz^{\alpha} \|_{L^\infty[0,1]}.
$$
However, $ r_{\alpha,k}(\zz)$ could be precomputed  and used without knowing the largest eigenvalue of $\calAt$.
Thus it eliminates the parameter $\delta$. In the applications $\delta$ is very small, but as shown 
in \cite{HLMM-19}, even when $\delta \approx 10^{-8}$ this may be beneficial to the Remez algorithm.
Since on $  [\delta,1]$ the function $z^\alpha$ has continuous derivatives, though getting very large 
at the left bound of the interval, Remez algorithm is less sensitive to round-off errors.
\end{remark}

The problem \eqref{bura}
has been studied extensively in the past, e.g. \cite{saff1992asymptotic,Stahl93,varga1992some}.  
Denoting the error by  
\beq\label{BURA-error}
E_{\alpha,k}:=\|r_{\alpha,k} (\zz) - \zz^{\alpha} \|_{L^\infty[0,1]},
\eeq
and applying \cite[Theorem 1]{Stahl93} we conclude that there is a constant 
$C_\alpha>0$,  independent of $k$, such that 
\begin{equation}\label{error-bound}
E_{\alpha,k}  \le C_\alpha e^{-2 \pi \sqrt{k \alpha}}.
\end{equation}
Thus, the BURA error converges exponentially to zero as $k$ becomes large.

Now the function  
$ u_{h,k}  \in V_h$ (and correspondingly its vector representation $ \tilwh \in \RR^N$) 
obtained from 
\beq
u_{h,k} 
=\lambda_{1,h}^{-\alpha} r_{\alpha,k} (\lambda_{1,h} \calAt^{-1}) f_h
\quad \mbox{or}  \quad 
\tilwh 
=\lambda_{1,h}^{-\alpha} r_{\alpha,k} (\lambda_{1,h} \wcalAt^{-1}) \tilfh
\label{wh}
\eeq
is called {\it fully discrete} approximation of    \eqref{eq:sol-function}.
Here $\calAt$ and $f_h$ are as in \eqref{fea} or \eqref{lumped-semi}
and $\wcalAt$ and $\tilfh$ are as in \eqref{fda}.

We stress that one does not need to know the exact value of $\lambda_{1,h}$. In fact, for any
$\delta >0$ such that $\delta  \le \lambda_{1,h}$,  the fully discrete solution 
\beq\label{delta-scale}
 u_{h,k}  =r_{\alpha,k} (\delta \calAt^{-1}) \delta^{-\alpha}  f_h
\eeq
represents  another good approximation to our problem.

In \cite{harizanov2020JCP}, we studied the error of these fully discrete solutions.
For the finite element case we obtain the error estimate 
\beq
\|u_{h} - u_{h,k}\| \le \lambda_{1,h}^{-\alpha} E_{\alpha,k} \|f_h\| \le C_\alpha \lambda_{1,h}^{-\alpha} e^{-2 \pi \sqrt{k \alpha}} \|f_h\| .
\label{pbest}
\eeq
with $\|\cdot\|$ the $L^2(\Omega)$-norm,
while in the finite difference case  we got
\beq
\|\tiluh - \tilwh  \|_{\ell_2} \le \lambda_{1,h}^{-\alpha} E_{\alpha,k} \|\tilfh \|_{\ell_2}  
\le C_\alpha \lambda_{1,h}^{-\alpha} e^{-2 \pi \sqrt{k \alpha}}\|\tilfh \|_{\ell_2} ,
\label{pbest-fd}
\eeq
where $\| \cdot \|_{\ell_2} $   
denotes the  Euclidean norm in $\RR^N$.
 
As an illustration, in Table \ref{t:error} we provide the computed  error $E_{\alpha,k}$ 
for some particular values of $\alpha$ and various $k$.
It is remarkable that for $\alpha =0.75$ one can get an error of the order $10^{-6}$ just for $k=6$.
However, for small values of $\alpha$ one needs high order rational functions to get a
reasonable error. Then by \eqref{pbest} and \eqref{pbest-fd} one gets a  
bound of the fully discrete error.

\renewcommand{\arraystretch}{1.2}
\begin{table}[ht!]
\caption{ \footnotesize The error $E_{\alpha,k}$ of BURA of $\zz^\alpha$, 
$ \zz \in[0,1]$}
\label{t:error}
\centering
\begin{tabular}{|c|c|c|c|c|c|c|}
\hline
$\alpha$&$E_{\alpha,5}$&$E_{\alpha,6}$&$E_{\alpha,7}$&$E_{\alpha,8}$&$E_{\alpha,9}$&$E_{\alpha,10}$\\ \hline
  0.75 & 2.8676e-5 & 9.2522e-6 & 3.2566e-6 & 1.2288e-6 & 4.9096e-7 & 2.0584e-7\\
  0.50 & 2.6896e-4 & 1.0747e-4 & 4.6037e-5 & 2.0852e-5 & 9.8893e-6 & 4.8760e-6\\
  0.25 & 2.7348e-3 & 1.4312e-3 & 7.8650e-4 & 4.4950e-4 & 2.6536e-4 & 1.6100e-4\\\hline
\end{tabular}
\end{table}

\subsection{The BURA solution method}
\label{sec:6.2}
Now we need to show that after finding  $r_{\alpha,k}(\zz)$ we can efficiently implement 
the computations of the solution by \eqref{wh}. This is possible due to the useful properties
of $r_{\alpha,k}(\zz)$, which could be found e.g in \cite{saff1992asymptotic,Stahl2003}.

 It is known that the best rational approximation
$r_{\alpha,k}(\zz)=P_k(t)/Q_k(\zz)$  of $\zz^\alpha$  for $\alpha \in (0,1)$ is
non-degenerate, i.e., the polynomials $P_k(\zz)$ and $Q_k(\zz)$ are of full
degree $k$.   Denote  the roots of $P_k(\zz)$ and $Q_k(\zz)$ by
$\zeta_1, \dots, \zeta_k$ and  $d_1, \dots, d_k$, respectively.    
It is shown in
\cite{saff1992asymptotic,Stahl2003}  that the roots are real, interlace and satisfy
\begin{equation}\label{interlacing}
0 > \zeta_1 > d_1 > \zeta_2 > d_2 > \cdots > \zeta_k > d_k.
\end{equation}
We then have
\beq
r_{\alpha,k} (\zz)=b\prod_{i=1}^k\frac{\zz-\zeta_i}{\zz-d_i}\label{prodr}
\eeq
where, by \eqref{interlacing} and the fact that
$r_{\alpha,k}(\zz)$ is a best approximation to a non-negative function, $b>0$ and
$P_k(\zz)>0$ and $Q_k(\zz)>0$ for $ \zz \ge 0$. 

Knowing the poles $d_i$, $i=1, \dots, k$ we can give an equivalent representation of \eqref{prodr}
as a sum of partial fractions, namely
\beq
  r_{\alpha,k}(\zz)=c_0+\sum_{i=1}^k \frac{c_i}{\zz- d_i}
  \label{trg}
  \eeq
  where $c_0>0$ and $c_i<0$ for $i=1,\ldots,k$.

Now changing the variable $\tt =1/\zz$ in $r_{\alpha,k}(\zz)$ we get a rational function 
$\trg(\tt)$ defined by
\beq\label{eq:error-tilde}
\trg(\tt):=r_{\alpha,k}(1/\zz)=\frac {\widetilde P_k(\tt)}{\widetilde Q_k(\tt)}.
\eeq
Here $\widetilde P_k(\tt)=\zz^k P_k(\zz^{-1})$ and $\widetilde Q_k(\tt)= t^k Q_k(\zz^{-1})$  and hence their coefficients 
are defined by reversing
the order of the coefficients in $P_k$ and $Q_k$ appearing in $r_{\alpha,k}(\zz)$.  
In  addition, \eqref{interlacing} implies that we have the following properties for  
the roots of $\widetilde P_k$ and $\widetilde Q_k$, 
$\widetilde d_i=1/d_i$ and $\widetilde \zeta_i = 1/\zeta_i$,  respectively.
\begin{equation}\label{interlacing1}
  0 >  \widetilde d_k > \widetilde \zeta_k > \widetilde d_{k-1} > \widetilde
  \zeta_{k-1}\cdots >   \widetilde d_1>\widetilde \zeta_1.
\end{equation}  
As a result we have the following lemma, cf. \cite{harizanov2020JCP}:
\begin{lemma}\label{lemma:BURA}  Let 
  $\widetilde c_0 = c_0 - \sum_{i=1}^k c_i \widetilde d_i=r_{\alpha,k}(0)=E_{\alpha,k}>0$, 
   and 
  $\widetilde c_i = - c_i { d}_i^{-2}  = - c_i {\widetilde d}_i^{2} >0$,  $i=1,\ldots,k$. Then for  $\alpha\in (0,1)$, 
\beq
\label{r-compute}
\trg(\tt)={\widetilde c}_0+\sum_{i=1}^k {\widetilde c}_i (\tt-{\widetilde d}_i)^{-1}.
  \eeq
  
\end{lemma}

From \eqref{eq:error-tilde} and \eqref{r-compute} we see that with $\xi$ replaced by $ \lambda_{1,h}^{-1} \calAt$ we have 
$$
r_{\alpha,k}(\lambda_{1,h} \calAt^{-1})= \widetilde r_{\alpha,k}(\lambda_{1,h}^{-1} \calAt)=
{\widetilde c}_0 \calIt +\sum_{i=1}^k {\widetilde c}_i (\lambda_{1,h}^{-1} \calAt-{\widetilde d}_i)^{-1}.
$$
Thus, the solution \eqref{wh} could be expressed by
$$
u_{h,k} = {\widetilde c}_0  f_h  + \sum_{i=1}^k {\widetilde c}_i w_i , \ \   
$$
where 
$w_i$ are solutions of the following $k$ systems 
$$
( \calAt-{\lambda_{1,h}}{\widetilde d}_i)w_i=  \lambda_{1,h}^{1-\alpha} f_h, \ \ i=1,\dots,k  .
$$
Note that ${\widetilde d}_i <0$ so that the corresponding matrix is 
positive definite and the summation is stable since ${\widetilde c}_i$ are all positive.

\subsection{BURA and URA methods for fractional diffusion reaction problems}
\label{sec:6.3}
The methodology, developed in Section~\ref{sec:6.2} can be straightforwardly extended to a generalization of \eqref{fda}, namely
\begin{equation}\label{DifReact}
({\wcalAt}^\alpha +\qq\II) \bfu =   \bff,\qquad \qq \ge 0,\quad \alpha\in(0,1).
\end{equation}
Such a problem appears for example in finite difference discretization of sub-diffusion-reaction elliptic problems or 
transient sub-diffusion problems. In the first case, $q$ corresponds to the reaction term. 
In the second case it is inversely proportional to the time step discretization, e.g., $\qq=\tau^{-1}$
if the backward Euler discretization in time is applied. 
Obviously, the solution can be written as $ \bfu = (\calI + \qq \calA^{-\alpha}) \calA^{-\alpha} \bff $.

Now for $\qq>0$ we introduce $\rqa(z)$ as the best uniform rational approximation of the function 
$g_\qq(z;\alpha):=\frac{z^\alpha}{1+\qq z^\alpha}$,   $z\in[0,1]$, that is 
\begin{equation}\label{ura}
 \rqa(\zz) := \argmin_{s(\zz)\in \rat k}\,  
 \max_{\zz \in [0,1]} \left |  s(\zz) - \frac{\zz^{\alpha}}{1+\qq \zz^\alpha} \right | .
\end{equation}

There is strong numerical evidence (see \cite{harizanov2019cmwa}) that $\rqa(z)$ inherits all useful properties of $r_{\alpha,k}(z)$, 
so the relations \eqref{interlacing}--\eqref{r-compute} remain valid. Furthermore, the corresponding version of 
the error bound \eqref{pbest-fd} reads as
\begin{equation}\label{Er-q}
\|\bfu- \bfu_k\|_{\ell^2}\le \lambda_{1,h}^{-\alpha}\Eqa\|\bff\|_{\ell^2}, 
\end{equation}
with $  \Eqa:=\|\rqa(z)-g_\qq(z;\alpha)\|_{L^\infty[0,1]}$.

The estimate of $\Eqa$ is obtained from the elements $\brqa(z)$ of a uniform rational 
approximation (URA) of $g_\qq(z;\alpha)$  defined as:
$$
\brqa(z)=g_{\qq_2}\left(r_{\qq_1,\alpha,k}(z),1\right)=\frac{r_{\qq_1,\alpha.k}(z)}{1+\qq_2r_{\qq_1,\alpha,k}(z)},\quad \qq=\qq_1+\qq_2,\; \qq_1,\qq_2\ge 0.
$$
Note that for all choices of $\qq_1,\qq_2\ge 0$, $\brqa(z)$ are rational functions in $\rat k$. 
However, these are NOT BURA-elements, unless $\qq_2=0$. Nevertheless, they approximate $g_\qq(z;\alpha)$ well
and the following estimate is valid, see \cite[Theorem 2.4]{harizanov2019cmwa}:
\begin{equation}\label{error_bounds}
\frac{(1+\qq_1)^2 
}{(1+\qq)^2}\le \frac{\Eqa}{E_{\qq_1,\alpha,k}}  \le \frac{1}
{1+\qq_2 E_{\qq_1,\alpha,k}}
, \;\forall \qq_1,\qq_2\ge 0,\; \qq_1+\qq_2=\qq.
\end{equation}
It follows that $E_{\qq,\alpha,k} < E_{\qq_1,\alpha,k}$ for $\qq_1<q$, so that 
$\Eqa$ monotonically decreases as $\qq$ increases. 
Various numerical experiments also show that $(1+\qq)\Eqa$ monotonically increases as $q$ increases for all values  
of the parameters $k,\alpha$, and $\qq$. In this case, $\lim_{\qq\to+\infty}(1+\qq)\Eqa=C(k,\alpha) E_{\alpha,k}$, 
with $1<C(k, \alpha) \le E_{\alpha,k}^{-1}$. This indicates that $ C(k, \alpha)$ could grow as $k$ grows.  
In practice, as could be seen from Table \ref{tab63},
$C(k,\alpha)$ changes linearly for realistic values of $k$.
\begin{table}
\caption{Numerical values of $C(k,\alpha) \approx 
(1+\qq)\Eqa / E_{\alpha,k}$ for  $\qq=400$ and various $\alpha$ and $k$}\label{tab63}
\begin{tabular}{|c|c|c|c|}
\hline
\multirow{2}{*}{$k$}&\multicolumn{3}{|c|}{$401 E_{400,\alpha,k}/E_{0,\alpha,k}$}\\ \cline{2-4}
& $\alpha=0.25$ & $\alpha=0.50$ & $\alpha=0.75$\\\hline
$3$ & 2.931 & 4.785  & 6.368\\\hline
$4$ & 3.774 & 6.588  & 8.967 \\\hline
$5$ & 4.708 & 8.614  & 11.813\\\hline
$6$ & 5.731 & 10.807 & 14.804\\\hline
\end{tabular}

\end{table}

Computing URA $\rqa$ becomes numerically unstable as $\qq\to\infty$, especially for small values of $\alpha$. 
This is due to the clustering at zero of the extreme points of the error function $\rqa(z)-g_q(z;\alpha)$, leading to necessity 
of execution of the Remez algorithm with higher than double-quadruple arithmetic precision. On the other hand, the set of 
extreme points of $\brqa(z)-g_\qq(z;\alpha)$ coincide with those of $r_{\qq_1,\alpha,k}(z)-g_\qq(z;\alpha)$ and the corresponding 
partial fraction decomposition \eqref{trg} of $\brqa$ can be cheaply derived from the one of $r_{\qq_1,\alpha,k}$. 
Thus, although not optimal, an URA element can be a useful approximation tool in practice.

Moreover, the error function in the URA case is not equi-oscillating. Its largest absolute value is at zero 
(the first extreme point) and the absolute value monotonically decreases with every successive extreme point. 
As a result, in vicinity of $1$, all URA elements give rise to smaller errors than the BURA element. 
Furthermore, we have the relations
$$
\Eqa>\bEqa,\qquad \lambda_{N,h}^\alpha<\sqrt{(1+\qq)(1+\qq_1)}+1.
$$
Thus, if one uses an upper bound $ \bar \lambda$
for $\lambda_{N,h}$, namely, $\bar \lambda \ge \lambda_{N,h}$, then
whenever $\qq_1 > {\bar \lambda}^\alpha -2$, and $\qq>\qq_1$ the URA method, related to $\brqa$ 
will have a smaller  error that the BURA method, related to $\rqa$. 
This means that in practice there is no need for the latter to be computed.

%%%%%%%%%%%%%%%%%%%%%%%%%%%%%%%%%%%%%%%%%%%%%%%%%
\section{Computational efficiency: a comparative analysis}
\label{sec:7}

\setcounter{section}{7}
\setcounter{equation}{0}\setcounter{theorem}{0}

\subsection{Computational complexity}
\label{sec:7.1}
From computational point of view, the basic idea behind the surveyed  methods is to 
approximate the solution of the non-local fractional diffusion problem through 
systems with sparse symmetric and positive definite matrices. 
Assuming that some solver of optimal complexity is used for the sparse systems, 
the total computational complexity of the method is determined by the number of 
these systems.

To understand better the computational results we state the following error estimates 
for the fully discrete scheme for the lumped mass finite element approximation in space 
with continuous piece-wise linear functions over a uniform mesh established in 
\cite{harizanov2020JCP}.

\begin{theorem} \label{masslumped} 
(\cite[Corollary 4.3]{harizanov2020JCP}) Let $\Omega\subset \R^2$
and  suppose that the operator $\calT$ defined in Subsection \ref{ss:problem}
provides full regularity lifting, i.e. $ \|\calT f \|_{H^2} \le c \|f\|$.
Then for $f\in H^{1+\gamma}(\Omega)$, $\gamma>0$, the exact solution  $u = \calA^{-\alpha}f$ and   
the fully discrete solution 
$u_{h,k}=r_{\alpha,k}(\lambda_{1,h} \calAt^{-1}) \lambda_{1,h}^{-\alpha} f_h$  
(for $\calAt$ obtained from using mass lumping \eqref{A-lumped})
satisfy 
\begin{equation}\label{error-total}
\|u - u_{h,k} \|  \le %\|u - u_h \| + \|u_h - w_h \| \le 
C (h^{2\alpha}+h^{1+\gamma}) \|f\|_{H^{1+\gamma}(\Omega)} + \lambda_{1,h}^{-\alpha} E_{\alpha,k} \|f_h\|.
\end{equation}
\end{theorem}

As we see, the first part of the error comes from the finite element approximation 
of the problem and the mass lumping. The second part is the error due to the use 
of approximation method to solve the system with fractional power of the related matrix.

\begin{table}[ht!]
\centering
\caption{Relative errors of BURA and $\kk$-Q methods.}
\label{tab:2D-results-CB}
\renewcommand{\arraystretch}{1.2}
\begin{tabular}{|c|c|cc|cc|cc|}
\hline
$\alpha$& $h$ & \multicolumn{2}{|c|}{BURA $k=9$} & \multicolumn{2}{|c|}{$\kk$-Q (9 solves)} 
     & \multicolumn{2}{|c|}{$\kk$-Q ($\kk=\frac13$)}\\
& &  $\ell_2$ & $\ell_\infty$ &  $\ell_2$ & $\ell_\infty$ &  $\ell_2$ & $\ell_\infty$ \\ \hline
$0.25$ & 
  $2^{-9}$  & 2.292e-4 & 1.875e-3 & 1.040e-2 & 1.207e-2 & 1.371e-4 & 1.847e-3 \\
& $2^{-10}$ & 2.029e-4 & 1.339e-3 & 1.039e-2 & 1.152e-2 & 6.815e-5 & 1.305e-3 \\
& $2^{-11}$ & 1.939e-4 & 8.219e-4 & 1.038e-2 & 1.097e-2 & 3.388e-5 & 9.196e-4 \\
& $2^{-12}$ & 1.922e-4 & 7.451e-4 & 1.038e-2 & 1.069e-2 & 1.671e-5 & 6.413e-4 \\\hline
$0.50$ & 
  $2^{-9}$  & 1.013e-5 & 8.787e-5 & 2.835e-3 & 2.904e-3 & 8.058e-6 & 9.110e-5 \\
& $2^{-10}$ & 8.304e-6 & 4.742e-5 & 2.830e-3 & 2.902e-3 & 2.840e-6 & 4.559e-5 \\
& $2^{-11}$ & 8.263e-6 & 2.433e-5 & 2.829e-3 & 2.902e-3 & 1.033e-6 & 2.280e-5 \\
& $2^{-12}$ & 8.291e-6 & 1.909e-5 & 2.828e-3 & 2.902e-3 & 4.118e-7 & 1.132e-5 \\\hline  
$0.75$ & 
  $2^{-9}$  & 6.110e-7 & 3.110e-6 & 1.502e-3 & 1.824e-3 & 7.118e-7 & 3.263e-6 \\
& $2^{-10}$ & 1.884e-7 & 1.037e-6 & 1.501e-3 & 1.823e-3 & 2.355e-7 & 1.198e-6 \\
& $2^{-11}$ & 1.500e-7 & 6.592e-7 & 1.500e-3 & 1.823e-3 & 1.138e-7 & 4.677e-7 \\
& $2^{-12}$ & 1.547e-7 & 4.574e-7 & 1.499e-3 & 1.823e-3 & 8.334e-8 & 2.079e-7 \\\hline 
\end{tabular}
\end{table}

The approach based on the number of sparse linear systems solves is used to compare the efficiency of  
BURA method (see \cite{HLMMV18,harizanov2019cmwa,harizanov2020JCP}) with the
method proposed by Bonito and Pasciak in \cite{BP-15} which is referred also as $k^\prime$-Q method.
The data in Table \ref{tab:2D-results-CB} are extracted from Table 2 of 
\cite{BP-15}. 

We consider a two-dimensional test problem with a Checker Board right-hand-side with
reference solution (taken as an exact solution) computed via FFT on a uniform square mesh with $h=2^{-15}$.
The first two sets of data concern the BURA (as defined in Section \ref{sec:6}) and  
$\kk$-Q method, both using $9$ linear system solves, while the  $\kk$-Q method with 
$\kk=1/3$, uses $120$ linear system solves for $\alpha=0.25,0.75$ and  
$91$ system solves for $\alpha=0.5$. 
Here we report the relative errors in $\ell_2$ and $\ell_\infty$ norms, namely
$$
\| e \|_{\ell_2}=   \|u - u_{h,k} \|_{\ell_2}/  \|  f_h\|_{\ell_2} \quad \mbox{and} \quad \| e \|_{\ell_\infty} =\|u - u_{h,k} \|_{\ell_\infty}/ \|  f_h\|_{\ell_\infty}, 
$$
where $ \| e \|_{\ell_2} $ is the standard Euclidean norm of the vector obtained from sampling $  e(x) =u(x) - u_{h,k}(x) $
at the mesh points and $ \| e \|_{\ell_\infty}$ is the maximum value of $e(x)$ at the mesh.

We see that for $\alpha=0.25$ and $\alpha=0.5$ and $9$ system solves (equivalent to $k=9$) the error is
essentially due to the rational approximation, it does not change when decreasing the mesh-size
(these are columns 3-6 in the table). For the $\kk$-Q-method with 120 solves
we see that the error decreases 
when decreasing the mesh-size $h$. This indicates that the finite element method error 
dominates.
However, the BURA error  is almost 50 times smaller than the error of the 
$\kk$-Q-method with the same number  of linear system solves.

Even in the case of worst approximation, 
$\alpha = 0.25$, BURA produces a reasonable error in the range of 
$10^{-4}$ when using only 9 system solves. Moreover, for a mesh-size 
$h = 2^{-9}$ BURA is outperforming Q-method on all meshes for $\alpha = 0.5$ 
and $\alpha = 0.75$. In contrast, the $\kk$-Q-method gives the same accuracy, 
but needs 91 and 120 system solves, respectively.

\begin{figure}[ht!]
\includegraphics[width=0.75\textwidth]{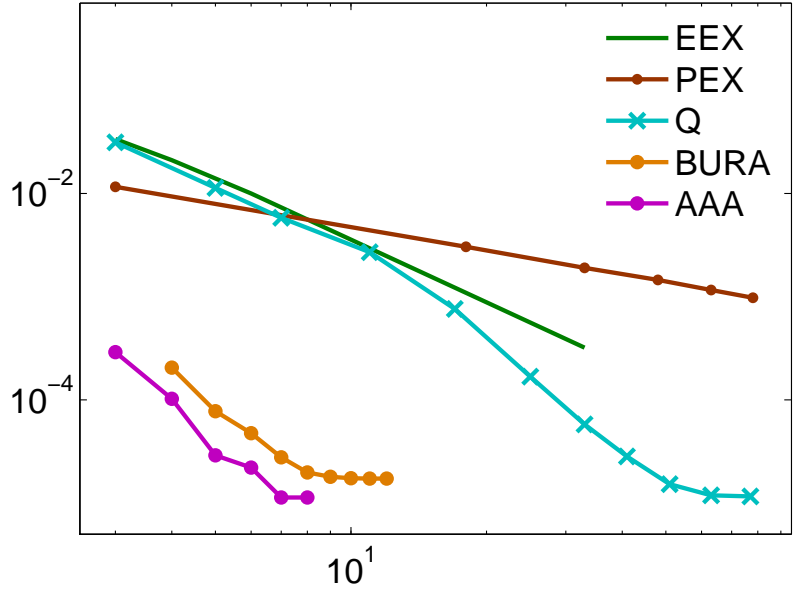}
\caption{Accuracy versus degree of rational approximation $k$ (number of linear systems solves) for
$\alpha=0.5$ and $h=1/512$.}\label{unicomp}
\end{figure}

Recently, a unified view  to the methods discussed in 
Sections \ref{sec:4} - \ref{sec:6} was presented in \cite{HOFREITHER2019}\footnotemark.
%{$^{(1)}$}
\footnotetext{ In this paper,  when talking about BURA, the author has had in mind 
the first variant of the method from \cite{HLMMV18} which is not robust with respect 
to the condition number of the matrix $\wcalAt$.}
The work is based on the 
observation that each of discussed above methods can be interpreted as generating  some 
rational approximations of $\wcalAt^{-\alpha}$ in the form
\begin{equation}\label{unirat}
\wcalAt^{-\alpha} \approx 
{\tilde c}_0 {\mathbb I} +\sum_{i=1}^k {\tilde c}_i (\wcalAt-{\tilde d}_i 
{\mathbb I})^{-1},
\end{equation}
where ${\tilde c}_i  \ge 0$ and ${\tilde d}_i<0$.
Thus, based on (\ref{unirat}), one can easily compare the efficiency of 
all methods considered in this survey. Such a comparison is provided
in Figure \ref{unicomp} (most of the data is from \cite{HOFREITHER2019})  
for $\alpha=0.5$, where the accuracy versus degree
of the rational approximation $k$ is displayed. Here we consider the test problem 
$(-d^2/dx^2)^\alpha u(x) =1$ for $x\in(-1,1)$ with boundary conditions $u(-1)=u(1)=0$.  
The discretization is done by linear finite elements with mass lumping 
(equivalent to a three-point finite difference approximation) on a uniform mesh with mesh-size $h=1/512$. 
The experiments are representative in the sense that the error estimates 
are independent of the space dimension $d$.
On Figure \ref{unicomp} we put together  the numerical results of BURA from \cite{harizanov2019cmwa,harizanov2020JCP}, 
$\kk$-Q-method from \cite{BP-15}, and the results of three additional methods labeled as:
EEX for the elliptic extension method 
\cite{BBNOS-18} with linear FEM discretization in the auxiliary direction $y$;
PEX for the pseudo-parabolic extension  \cite{Vabishchevich15} with Crank-Nicolson 
discretization in time  $t$; AAA, based on the Chebfun generated best uniform 
rational approximation of a discrete set of values of $z^\alpha$ on 
$(\lambda_{1,h}, \lambda_{N,h})$ by using  the Symbolic Math Toolbox of
MATLAB R2008b, \cite{mathworks09}, introduced and discussed in
\cite{HOFREITHER2019}. The data for these three methods are extracted from
Figure 2 of \cite{HOFREITHER2019}. In agreement with the theory,
for any fixed $k$, BURA outperforms EEX, PEX and $\kk$-Q methods.
We note that rather small $k$ is sufficient to balance the discretization error of order $O(h^{2\alpha})$
(the case if  $ f \in L^2(\Omega)$) with BURA error.

As discussed in Section \ref{sec:6.1}, the computation of BURA for 
larger $k$ faces certain problems due to the computational instability 
of the Remez algorithm. The AAA method avoids this difficulty under
the assumption that the extremal eigenvalues of $\wcalAt$ are available.
The so-called “adaptive Antoulas-Anderson” (AAA) \cite{NST-18} algorithm 
exploits representation of the rational approximant in barycentric form and
greedy selection of the support points. The method proposed in  
\cite{HOFREITHER2019} is based on AAA approximation of $z^{-\alpha}$ for
$z \in [\lambda_{1,h},\lambda_{N,h}]$, while the BURA method essentially uses 
the approximation on $[\lambda_{1,h}, \infty)$. According to Remark \ref{delta-BURA},
the error of AAA method is always smaller than the error of BURA method.
This is also in agreement with the computations of the BURA by Remez algorithm 
presented in \cite{HLMM-19}. Figure \ref{unicomp} shows that the errors of 
AAA and BURA have similar behavior with respect to $k$. 
When $\lambda_{1,h}$ and 
$\lambda_{N,h}$ are known, the AAA approximation is cheaper to get
for different $\alpha$ and the computations are more stable for larger $k$.
Moreover, if we have bounds $ \underline \lambda$ and $ \overline \lambda$, 
such that  $\underline \lambda \le \lambda_{1,h}$ 
and $\lambda_{N,h} \le \overline \lambda$, then the application of AAA method 
will still generate a good approximation to $z^{-\alpha}$ on $(\underline \lambda, \overline \lambda)$,
which  gives the AAA method some practical advantage.

\subsection{Parallel efficiency}
\label{sec:7.2}
The development of highly efficient parallel algorithms for large-scale problems is
a topic of rapidly growing interest. In the case of fractional diffusion problems, 
the parallel implementation is of even stronger importance. 

The following example illustrates how extreme computational demands could appear. 
Let the problem in $\Omega\subset\RR^3$ be discretized by 
linear finite elements with a mesh parameter $h$, and let $f(x)\in L^2(\Omega)$. 
Then: (i) Standard diffusion: FEM error estimate: $O(h^2)$; $h=10^{-2}$ is needed for
$O(10^{-4})$ accuracy, leading to $N = O(10^6)$;
(ii) Modest fractional diffusion: FEM error estimate: $O(h^{2\alpha})$; for
$\alpha = 0.5$, $h = 10^{-4}$ is required to get $O(10^{-4})$ accuracy, leading to
$N = O(10^{12})$. The last size is a challenge even for the modern supercomputers.

The first study on parallel solution of fractional diffusion problems was published
in \cite{KMV-18}. The fractional Laplacian in the unit cube $\Omega= (0,1)^3$ is
considered, where the seven point stencil is applied to discretize the Laplace operator.
The $\kk$-Q-method with $\kk=1/3$ is used which means 91 auxiliary systems
with sparse symmetric and positive definite matrices. The PCG is utilized as a basic 
iterative solution method for these systems where a parallel multi-grid (MG) implementation 
from the Trilinos ML package is the preconditioner. 
\begin{table}[ht!]
\caption{Parallel scalability: maximum 61 physical cpu cores per node.}
\label{ParE}
\begin{center}
\renewcommand{\arraystretch}{1.2}
\begin{tabular}{|c|rr@{~~~~}rr|r@{~~~~}r|}
\hline 
 & & \multicolumn{2}{c}{$N=128^3$} & ~~~~~~~~ & \multicolumn{2}{c|}{$N=256^3$}\\ 
\hline
  Nodes & & T[s] & ~~~E[\%] & & T[s] & ~~~E[\%] ~~~ \\
\hline
  1 & & 146 & & & 989 & \\
  2 & & 59 & 124 & & 455 & 109 \\
  4 & & 37 & 98 & & 244 & 101 \\
\hline
\end{tabular}
\end{center}
\end{table}
The reported parallel times T[s] and efficiencies E[\%] for $\alpha=0.5$ are shown in 
Table \ref{ParE}. We stress that in this implementation, the distribution of 
the 91 solves between the nodes is optimized, taking into account the different number of PCG iterations
for each of them, needed to reach the stopping criteria of $10^{-10}$.

Various  aspects of the parallel 
implementation of the surveyed methods are discussed in \cite{CSMK-18,CSMK-19,MRAplus-19}\footnotemark. 
\footnotetext{ In these papers, the authors have used some of the earlier variants of BURA 
from \cite{HLMMV18,HLMMP-19}, which are not robust with respect to the condition number 
of the matrix $\wcalAt$.}  
A scalability analysis of the PEX, $\kk$-Q and BURA methods is presented in 
\cite{CSMK-19}, where the test problem in $\Omega=(0,1)^3$ with a CheckerBoard 
right-hand-side is considered with up to $512^3$ unknowns. 
The EEX method is excluded from the list of studied methods, as less suitable
due to the high memory requirements in the 3D case.
The discussed question is which parallel algorithm is recommended to achieve a certain accuracy 
for a given $\alpha\in\{0.25,0.75\}$. 

A less commonly used approach to a comparison analysis of parallel efficiency
of EEX, PEX, $\kk$-Q and BURA methods is proposed in \cite{MRAplus-19}. The presented 
results are based on using up to 32 nodes of the supercomputer, with a setting 
of up to 16 cpu cores per node\footnotemark. 
For example, in \cite[Table 6.8]{MRAplus-19}, the best parallel times to 
achieve a given accuracy are shown. Similarly, \cite[Table 6.9]{MRAplus-19} displays the best speed-ups 
versus accuracy. The performed analysis has shown that the selection of the best
algorithm is problem dependent.

\footnotetext{ The parallel numerical tests discussed above are performed on the 
Bulgarian Academy of Sciences supercomputer Avitohol (http://www.iict.bas.bg/avitohol/). }

At the end of this section, it is worth to point out that all considered methods 
have been implemented in their original formulations. Now, after the  unified interpretation
of all of them (see \cite{HOFREITHER2019}) as certain rational approximations, the related additive
representation as a sum of partial fractions is expected to be used in the development of
future parallel algorithms.

%%%%%%%%%%%%%%%%%%%%%%%%%%%%%%%%%%%%%%%%%%%%%%%%%
\section{Challenges beyond the scalar elliptic case}
\label{sec:8}

\setcounter{section}{8}
\setcounter{equation}{0}\setcounter{theorem}{0}

\subsection{Time dependent space-fractional diffusion problems}
\label{sec:8.1}
Let us consider the time dependent problem: find $u(x,t)$ for 
$(x,t)\in\Omega\times (0,T]$ such that 
\beq\label{eq:timed-c}
\frac{\partial u(x,t)}{\partial t} + \calA^\alpha u(x,t) = f(x,t)
\quad \mbox{and} \quad u(x,0) = v(x), 
\eeq
with $v(x)$ and $T >0$ given initial data, and $T >0$ a real number.
In this section, for simplicity of the presentation we will assume
that the finite difference method is used to approximate $\calA$ 
in space. We introduce also a uniform mesh in time with step size
$\tau=T/M$, $M$ is a given integer parameter. Following the established
matrix notations we write the fully explicit two-level scheme in the
form:
\beq\label{eq:FE-2LS}
\frac{\bfu^{j+1}-\bfu^j}{\tau} + \wcalAt^\alpha\bfu^j = \mathbf {F}^j,
\quad j=1,\ldots, M,
\eeq
where the upper index $j$ indicates the related mesh
function (vector) at the time level $t=j\tau$. This scheme is conditionally 
stable. 

The
following regularized scheme is
proposed in \cite{V-16}
\beq\label{eq:RE-2LS}
({\mathbb I} + \tau {\mathbb R})  \frac{\bfu^{j+1}-\bfu^j}{\tau} + 
\wcalAt^\alpha\bfu^j = \mathbf {F}^j,
\quad j=1,\ldots, M,
\eeq
${\mathbb R}=\sigma(\alpha\wcalAt + (1-\alpha){\mathbb I})$,
proving unconditional stability if $\sigma\ge 0.5$.  To implement the 
scheme (\ref{eq:RE-2LS}) one has to perform matrix-vector multiplication
with $\wcalAt^\alpha$. For this purpose the representation 
$\wcalAt^\alpha = \wcalAt \wcalAt^{-(1-\alpha)}$ is used and reformulation 
of the fractional problem to a pseudo-parabolic (see Subsection \ref{sec:4.2})
is applied to approximate the solution of systems with $\wcalAt^{1-\alpha}$. 

The alternative approach proposed in \cite{AN-17} is based on rational approximation
of the discrete fractional Laplacian.
It is obtained from applying the Gauss-Jacobi quadrature to the integral representation
of $\wcalAt^\alpha$ in $[0,1]$ and is described in Subsection \ref{sec:5.2}.

The methods form \cite{AN-17,V-16} have promising stability properties that are
confirmed by numerical tests. The drawback for large scale problems
in space is that their accuracy is not robust with respect to the condition number
$\kappa(\wcalAt)$. The matrix-vector multiplication with $\wcalAt^\alpha$ could be  
avoided by applying the unconditionally stable backward Euler scheme 
\beq\label{eq:BE-2LS}
\frac{\bfu^{j+1}-\bfu^j}{\tau} + \wcalAt^\alpha\bfu^{j+1} = \mathbf {F}^j,
\quad j=1,\ldots, M.
\eeq
The implementation of (\ref{eq:BE-2LS}) requires solution of
linear systems with the matrix $\wcalAt^\alpha + \frac{1}{\tau} {\mathbb I}$ at 
each time step $j=1,\ldots, M$. For this purpose, one can use the BURA
method for fractional diffusion reaction problems discussed in
Subsection \ref{sec:6.3} with $\qq=\frac{1}{\tau}$. The BURA error estimate is
robust with respect to $\kappa(\wcalAt)$. In this setting, there appear new challenges
related to the numerical stability of the Remez algorithm for $\qq \gg  1$. 
Some combination of BURA and URA methods could be helpful in this context.

\subsection{Coupled problems involving fractional diffusion operators}
\label{sec:8.2}

A majority of the real-life applications are described by coupled 
problems. Among many others, we could mention the fractional diffusion 
epidemic models,  \cite{BDG-09}, the two-phase flow models based on the 
Navier-Stokes equations combined with a fractional Allen-Cahn 
mass-preserving model, \cite{SXK-17}, or surface quasi-geostrophic flows, \cite{BN-20,P-geophysical}. 
To illustrate some basic ideas and the related challenges 
we will consider the system of time dependent fractional-in-space 
diffusion-reaction equations for the unknown functions $u_\ell(x, t)$ 
in the form
\beq\label{eq:coups-c}
\frac{\partial u_\ell(x,t)}{\partial t} + \calA_\ell^{\alpha_\ell} u(x,t) = 
{\mathcal R}_\ell(u_1,\ldots,u_m) + f_\ell (x,t), \quad \ell=1,\ldots,m,
\eeq
 with given initial data $u_\ell(x,0) = u_{\ell,0}(x)$, $ \ell=1,\ldots,m$.
The system (\ref{eq:coups-c}) is coupled trough the reaction operators
${\mathcal R}_\ell(u_1,\ldots,u_m)$. Now, we rewrite the system in the
form of abstract Cauchy problem
\beq\label{eq:ACP-c}
\frac{\partial U(x,t)}{\partial t} = (\calA^\alpha+{\mathcal R})U(x,t) +F(x,t), 
\quad U(x,0)=U_0(x),
\eeq
where $\calA^\alpha:=diag(\calA_1^{\alpha_1}, \dots, \calA_m^{\alpha_m})$ 
and ${\mathcal R}$ are the fractional diffusion 
and the reaction, respectively, 
$U(x,t)=(u_1(x,t),\ldots,u_m(x,t))^T$,  $F(x,t) =(f_0(x,t), \dots, f_m(x,t))^T$ and 
$ U_0(x)=(u_{1,0}(x),\ldots,u_{m,0}(x))^T$.

Nowadays, the operator splitting is a commonly used approach in solving
such kind of problems. The basic ideas are associated
with the pioneering works of Yanenko, \cite{Yanenko-62}, Marchuk, 
\cite{Marchuk-68} and Strang, \cite{Strang-68}. In the case of
standard (not fractional) elliptic operator $\calA$,  i.e.
$\alpha_\ell \equiv 1$, the advantages of second (or higher) order splitting 
methods are well understood.  As a principle, they use Crank-Nicolson
like approximation of the derivative in time, thus involving in particular 
matrix-vector multiplication with the discrete diffusion operator. 
As was discussed in the previous subsection, the development
of robust method for multiplication with $\calAt^\alpha$ is
still a challenging problem. This is the main reason to restrict our 
consideration to the following sequential splitting algorithm:

For $ (j-1)\tau < t \le j\tau $,  $j=1,2,\ldots , M$,  %\vspace{-1mm}
$$ 
\begin{array}{ll}
\quad {\displaystyle\frac{\partial{U_1^j(x,t)}} {\partial t}} =  {\calA^\alpha} {U}_1^j(x,t), 
\quad  {U}_1^j(x,(j-1)\tau) = {U}_2^{j-1}(x,(j-1)\tau), & \vspace{3mm}\\
\quad {\displaystyle \frac{\partial{U}_2^j(x,t)} {\partial t}}  =  {\mathcal R} {U}_2^j(x,t) +F(x,t), 
\quad   {U}_2^j(x,(j-1)\tau) = {U}_1^j(x,j\tau), &
\end{array}
$$
%\vspace{-1mm}\hspace{-1mm}
where $\tau = T/M$ is the time step and
$U_2(x,0)=U_0(x)$. The function $U^j_{sp}(x,j\tau)=U^j_2(x,j\tau)$
is a {\it {sequential splitting approximate solution}} of  (\ref{eq:ACP-c}). 

Here we follow the abstract convergence analysis from \cite{Farago-05}, 
assuming that the operators $\calA$ and $\mathcal R$ are bounded 
with respect to $t\in[0,T]$ and 
the abstract Cauchy problem (\ref{eq:ACP-c}) is well posed. Then, the above sequential splitting 
is unconditionally stable and the splitting error 
is $O(\tau)$ \cite[Theorem 1]{Farago-05}. 
We will assume also that the backward Euler time-stepping scheme for the fractional 
diffusion sub-problems is combined with a properly chosen Runge-Kutta 
solution methods of the sub-problems associated with the reaction operator 
${\mathcal R}$, thus ensuring the targeted accuracy of $O(\tau)$. 

There are several different errors in the composite algorithm. Their
balancing is of a key importance. Now, for simplicity of the presentation,
we will assume that $\alpha_\ell=\alpha$, and that a uniform mesh with 
mesh parameter $h$ is used for approximation of the diffusion operator.
Under certain usual assumptions, the convergence rate of discretization 
in space of the fractional diffusion problems is $O(h^{2\alpha})$, see, e.g. 
\cite{bonito2019sinc}. The application of BURA method in the case of
backward Euler time-stepping was discussed in the previous subsection,
see the paragraph after (\ref{eq:BE-2LS}). From the numerical data 
presented in \cite{HLMM-19} we can deduct that for $\qq \gg 1$ the BURA
error behaves like $O(\frac{1}{q})$.  This follows from the estimate
\begin{equation}\label{Spl-err}
E_{\qq,\alpha,k}= C(k,\alpha)E_{\alpha,k}/(\qq+1) = O\left(
e^{-2\pi\sqrt{k\alpha}}/q \right ),
\end{equation}
which is concluded form \eqref{error_bounds} and the analysis there after 
(see also, \eqref{pbest} and recall that $E_{\alpha,k} =E_{0,\alpha,k}$), 
(\eqref{pbest-fd}) and Table \ref{tab63}). 
Thus, taking $\tau=1/\qq$ in (\ref{Spl-err}) and combining with (\ref{Er-q}), we get
the following asymptotic estimate holds true for the BURA 
error of the fractional diffusion sub-problem with backward 
Euler discretization in time
$O\left (\tau e^{-2\pi\sqrt{k\alpha}}\right )$, 
where $k$ is the degree of the best uniform rational approximation. Thus we get that
the considered composite sequential splitting algorithm has a total error decay $O(\tau + h^{2\alpha})$.

%%%%%%%%%%%%%%%%%%%%%%%%%%%%%%%%%%%%%%%%%%%%%%%%%
\section{Concluding remarks}
\label{sec:9}

In this survey we discussed various numerical methods for solving 
equations \eqref{semi} arising in discretization of  fractional  by powers of multidimensional elliptic problems.

Though quite different
in derivation and error analysis these methods have one common underlying feature: they all produce 
some rational approximation 
$ r_k(\calAt) = P_k(\calAt)/Q_k(\calAt)$ of $\calAt^{-\alpha}$ so that instead $  u_h = \calAt^{-\alpha} f_h  $ we compute 
$   u_{h,k} = r_k(\calAt) f_h $. Using spectral argument we see easily that  the error $ u_h - u_{h,k}$ is estimated 
by the error $   \max_{\zz \in [\lambda_{1,h}, \lambda_{N,h}]} |\zz^{-\alpha} - r_k(\zz)|$.
Thus, one concludes that any ``good" approximation of $\zz^{-\alpha}$ on $  [\lambda_{1,h}, \lambda_{N,h}]$
will produce a solution of \eqref{semi} as well.
This is equivalent to finding an approximation of $\zz^\alpha$ on 
$[\lambda_{N,h}^{-1}, \lambda_{1,h}^{-1}]$, which upon introducing a scaling  of $\calAt$ by $ \lambda_{1,h}^{-1}$,
is reduced to minimization in $[\epsilon, 1]$, with $\epsilon=\lambda_{1,h}/\lambda_{N,h}$. 

Since $\epsilon$ is very small
(it diminishes like $\min h^2$) we essentially need to find a ``good" approximation on $(0,1]$ to $\zz^\alpha$, which has 
singular derivative at $0$.   Remez algorithm for computing this approximation 
becomes more numerically unstable and computationally expensive for small
$\alpha$ and/or large $k$. Due to the theoretical results \cite[Theorem 3]{saff1992asymptotic},
both the zeros $\{\zeta_i\}_1^k$ and poles $\{d_i\}_1^k$ of the BURA $r_{\alpha,k}(\zz)$ of $\zz^\alpha$, $\zz\in[0,1]$, 
cluster at zero, when $k$ increases. More precisely, for every choice of $-\infty\le a\le b <0$ 
$$
\lim_{k\to\infty}\frac{1}{\sqrt{k}}card\{\zeta_i\in[a,b]\}=\frac{\sqrt{\alpha}}{\pi}\int_{|b|}^{|a|}\frac{dt}{t\sqrt{1+t}}
=\lim_{k\to\infty}\frac{1}{\sqrt{k}}card\{d_i\in[a,b]\}.
$$
In other words, the number of poles (as well as zeros) on any given interval $(a,b)$, $a<b<0$, grows like $\sqrt{k}$ and since 
the total number is $k$, this proves that for large 
$k$ an $O(\sqrt{k})$ of the poles (as well as the zeros) of $r_{\alpha,k}(\zz)$
are as close to the origin as one wishes. Similar result, \cite[Theorem 4]{saff1992asymptotic}
is valid for the extreme points $\{\eta_i\}_1^k$ of the error $r_{\alpha,k}(\zz)-\zz^\alpha$.
To illustrate the clustering, we give below the distribution of the poles $d_j$, $j=1, \dots,k$, of $r_{\alpha,k}(t)$ for $k=8$
(see, \cite[Table 38]{HLMM-19}): \\
%{}\\
$\alpha =0.25$:
$ d_1=-$2.39 e-11,  $d_2=-$8.37 e-9,  $d_3=-$5.95 e-7, $ \dots, d_8=-$1.53;\\$\alpha=0.50$:
$d_1=-$7.35 e-8,  $d_2=-$3.98 e-6,  $d_3=-$7.62 e-4,  $\dots, d_8=-5.43$;\\ 
$\alpha=0.75$:
$d_1= -$2.38 e-6,  $d_2=-$5.93 e-5, $d_3= -$6.50 e-4,  $\dots, d_8=-17.9$.

The clustering of the poles (the extremal points as well) shows that, 
high numerical accuracy and computer arithmetic precision is needed for 
computing $r_{\alpha,k}(\zz)$, when $\alpha \ll 1$ and/or $k\gg1$.
This is the most serious challenge for Remez algorithm which exhibits instability while  computing 
$r_{\alpha,k}(\zz)$ for $k \ge 10$. 

In Section \ref{sec:7.1} we discussed the AAA algorithm. It has been used to generate 
a best uniform rational approximation of a discrete set of values of $z^\alpha$, 
$z\in[\lambda_{1,h},\lambda_{N,h}]$. The Chebfun implementation of
AAA algorithm has been successfully used by Hofreither \cite{HOFREITHER2019}
for numerical solution of fractional diffusion problems. The reported numerical
results (see, e.g., Fig. \ref{unicomp}) show that this could be an attractive  
practical approach allowing to utilize the theoretical advantages of the BURA methods for
wide class of applications. But due to the clustering of the poles near zero even this method 
suffers of instabilities for large $k$. Based on this observation, Hofreithger in \cite{HOFREITHER2020}
presented a method for computing $r_{\alpha,k}(\zz)$ using radial basis functions and adaptation  
procedure. This resulted in an algorithm that is robust with respect to $\alpha$ and $k$
and is the best available method for fining approximately BURA in our opinion.

Moreover,  any of the methods discussed in this paper could be rated 
for accuracy and performance by checking whether they produce 
a rational approximation with clustering of the poles near the origin. For example, the 
method of Vabishchevich for solving a pseudo-parabolic equation  with constant 
time stepping (described in Subsection \ref{sec:4.2}) will be much less efficient than 
BURA method. An improvement that uses adaptive time-stepping and
produces certain clustering of the poles of the corresponding rational function is presented 
and justified in \cite{DLP-Pade}.

Finally, we conclude with some  
challenges. As noted in \cite{Tref-20}, the AAA rational approximation is 
near-best. By definition, it depends on the spectrum of the system’s matrix. 
In this sense, the method is not fully robust with respect to the condition number. 
For example, some further improvements are required to stabilize the  convergence, 
when applied to fractional diffusion reaction problem (\ref{DifReact}) for 
reaction coefficients $\qq \gg1$, a case  beyond the 
theoretical studies of Stahl \cite{Stahl2003}. The quality of the AAA approximation 
depends also from the data (the set of values of the related function) that are 
approximated. In \cite{HOFREITHER2019}, the values of $z^\alpha$ are taken on a 
uniform mesh of $[\lambda_{1,h},\lambda_{N,h}]$, in general, not 
the best choice.  One should expect some improvement if the exponential 
node clustering around the singularity point(s) is taken into account, see e.g., 
\cite{Tref-20}. The last, but not least comment concerns the 
Chebfun implementation. Although the computational stability of AAA algorithm is ensured, 
the opportunity to use double (and even higher) precision could be very useful 
in some more complicated (ill conditioned) problems.

Further progress in computational stability has been recently reported in 
\cite{HOFREITHER2020}. The proposed BRASI algorithm is based on the assumption that 
the best rational approximation of $f(z)$ must interpolate the function at a certain 
number of interpolation nodes $z_j$. It iteratively rescales the intervals $(z_{j-1},z_j)$ 
with the goal of equilibrating the local errors. The barycentric rational formula is 
used in the implementation. For example, results for $k=93$, $\alpha=0.25$ and $q=400$ 
are reported to illustrate the improved stability for $q \gg 1$.

\setcounter{section}{9}
\setcounter{equation}{0}\setcounter{theorem}{0}

%%%%%%%%%%%%%%%%%%%%%%%%%%%%%%%%%%%%%%%%%%%%%%%%%

\section*{Acknowledgments}
The partial support through the Grant No BG05M2OP001-1.001-0003, 
financed by the Science and Education for Smart Growth 
Operational Program (2014-2020) and co-financed by the EU 
through the European Structural and Investment Funds, and
through the Bulgarian NSF Grant No DFNI-DN12/1 is highly appreciated.
The work of  R. Lazarov has been partially supported by  NSF-DMS \#1620318 grant. \\

%%%%%%%%%% References %%%%%%%%%%%%%%%%%%%%%%%%%%%%%%%%

{\small

%\listoftables
% \listoffigures

}

%%%%%%%%%% put authors' addresses here, in \it %%%%%%%%


\begin{thebibliography}{99} %% {100}
\normalsize %%%%%%%%%%%%%%%%%%%%%%%

\bibitem{AN-17}
L. Aceto, P. Novati,
Rational approximation to the fractional Laplacian operator in
reaction-diffusion problems. {\em SIAM J. Sci. Comput.} \textbf{39},
No 1 (2017), A214--A228, DOI: https://doi.org/10.1137/16M1064714.

\bibitem{AN-18}
L. Aceto, P. Novati,
Efficient implementation of rational approximations to fractional
differential operators.
{\em Journal of Scientific Computing} \textbf{76}, No 1 (2018), 651--671, DOI: https://doi.org/10.1007/s10915-017-0633-2.

\bibitem{AN-19}
L. Aceto, P. Novati,
Rational approximations to fractional powers of self-adjoint
positive operators.
{\em Numerische Mathematik} \textbf{143} (2019), 1--16, DOI: https://doi.org/10.1007/s00211-019-01048-4.

\bibitem{AN-20}
L. Aceto, P. Novati,
Fast and accurate approximations to fractional powers of operators.
{\em arXiv:2004.09793} (2020).

\bibitem{Acosta-Borth2017}
G.~Acosta, J.P.~Borthagaray,
A fractional Laplace equation: regularity of solutions and
finite element approximations.
{\em SIAM J. Numer. Anal.} \textbf{55}, No 2 (2017), 472--495, DOI: https://doi.org/10.1137/15M1033952.

\bibitem{ABB-19}
G. Acosta, F.M. Bersetche, J.P. Borthagaray,
Finite element approximations for fractional evolution problems.
{\em Fract. Calc. Appl. Anal.}
\textbf{22}, No 3 (2019), 767--794; DOI: 10.1515/fca-2019-0042.

\bibitem{AM-20}
M. Ainsworth, Z. Mao,
Fractional phase-field crystal modelling: analysis,
approximation and pattern formation.
{\em IMA J. of Applied Mathematics} \textbf{85},
No 2 (2020), 231--262, DOI: 10.1093/imamat/hxaa004.

\bibitem{balakrishnan}
A. Balakrishnan, \newblock Fractional powers of closed operators
and the semigroups generated by them.
{\em Pacific J. Math.} \textbf{10}, No 2 (1960), 419--437.

\bibitem{BKM-19}
T. B{\ae}rland, M. Kuchta, K.-A. Mardal,
Multigrid methods for discrete fractional Sobolev spaces.
{\em SIAM J. Sci. Comput.} \textbf{41}, No 2 (2019), A948--A972, DOI: https://doi.org/10.1137/18M1191488.

\bibitem{BBNOS-18}
A. Bonito, J.P. Borthagaray, R.H. Nochetto, E. Ot{\'a}rola, A.J. Salgado,
Numerical methods for fractional diffusion. {\em Comput. Visual Sci.},
\textbf{19} (2019), 19--46, DOI: https://doi.org/10.1007/s00791-018-0289-y.

\bibitem{bonito2019numerical}
A. Bonito, W. Lei, J.E. Pasciak,
Numerical approximation of the integral fractional Laplacian.
{\em Numer. Math.} \textbf{142}, No 2 (2019), 235--278, DOI: https://doi.org/10.1007/s00211-019-01025-x.

\bibitem{BN-20}
A. Bonito, M. Nazarov,
Numerical simulations of surface-quasi geostrophic flows on periodic domains.
 {\em Preprint arXiv:2006.01180} (2020).

\bibitem{bonito2019sinc}
A. Bonito, W. Lei, J.E. Pasciak,
On sinc quadrature approximations of fractional powers of regularly
accretive operators.
{\em J. of Numerical  Mathematics} \textbf{27}, No 2 (2019), 57--68, DOI: 10.1515/jnma-2017-0116.

\bibitem{BP-15}
A. Bonito, J.E. Pasciak,
Numerical approximation of fractional powers of elliptic operators.
{\em Mathematics of Computation} \textbf{84}, No 295 (2015), 2083--2110,
DOI: https://doi.org/10.1090/S0025-5718-2015-02937-8.

\bibitem{BP-17}
A. Bonito, J.E. Pasciak,
Numerical approximation of fractional powers of regularly accretive
operators. {\em IMA J. Numerical Analysis} \textbf{37}, No 3 (2017), 1245--1273, 
DOI: https://doi.org/10.1093/imanum/drw042.

\bibitem{BDG-09}
D. Brockmann, V. David, A.M. Gallardo,
Human mobility and spatial disease dynamics.
{\em Reviews of Nonlinear Dynamics and Complexity} \textbf{2} (2009), 1--24.

\bibitem{caffarelli2007extend}
L. Caffarelli, L. Silvestre,
An extension problem related to the fractional Laplacian.
{\em Commun. in Partial Differential Equations} \textbf{32}, No 8 (2007), 1245--1260,
DOI: https://doi.org/10.1080/03605300600987306.

\bibitem{CaffareliStinga16}
L.A.~Caffarelli and P.R.~Stinga,  Fractional elliptic equations, {C}accioppoli estimates and regularity.
{\em Annales de l'Inst. Henri Poincare (C) Non Linear Analysis}, \textbf{33}, No 3 (2016), 767--807,
DOI: https://doi.org/10.1016/j.anihpc.2015.01.004.

\bibitem{CH-18}
A.L.~Chang, H.G.~Sun,
Time-space fractional derivative models for CO2 transport in
heterogeneous media. {\em Fract. Calc. Appl. Anal.} \textbf{21}, No 1 (2018), 151--173; 
DOI:  10.1515/fca-2018-0010.

\bibitem{CNOSalgado-2016}
L. Chen, R. Nochetto, O. Enrique, A.J. Salgado,
Multilevel methods for non-uniformly elliptic operators and fractional diffusion.
{\em Math. of Computation} \textbf{85} (2016), 2583--2607, DOI: https://doi.org/10.1090/mcom/3089.

\bibitem{ciarlet2002}
P.G. Ciarlet,
{\em The Finite Element Method for Elliptic Problems}.
Classics in Applied Mathematics, SIAM (2002).

\bibitem{CSMK-18}
R. \v{C}iegis, V. Starikovi\v{c}ius, S. Margenov, R. Kriauzien\'{e},
A comparison of accuracy and efficiency of parallel solvers for
fractional power diffusion problems.
In: {\em Parallel Processing and Applied Mathematics. PPAM 2017.
Lecture Notes in Computer Science} (Eds: R. Wyrzykowski, J. Dongarra, E. Deelman, K. Karczewski) 
\textbf{10777} (2018), 79--89.

\bibitem{CSMK-19}
R. \v{C}iegis, V. Starikovi\v{c}ius, S. Margenov, R. Kriauzien\'{e},
Scalability analysis of different parallel solvers for 3D fractional
power diffusion problems.
{\em Concurrency and Computation: Practice and Experience},
\textbf{31}, No 19 (2019), DOI: 10.1002/cpe.5163.

\bibitem{CV-19}
R. \v{C}iegis, P.N. Vabishchevich,
Two-level schemes of Cauchy problem method for solving fractional
powers of elliptic operators.
{\em Computers $\&$ Mathematics with Appl.} \textbf{80}, No 2 (2019), 305--315;
DOI: 10.1016/j.camwa.2019.08.012.

\bibitem{CV-20}
R. \v{C}iegis, P.N. Vabishchevich,
High order numerical schemes for solving fractional powers of elliptic
operators. {\em J. of Computational and Appl. Math.}
\textbf{372} (2020); DOI: 10.1016/j.cam.2019.112627.

\bibitem{Djrbashian:1993}
M.~M. Djrbashian,
{\em Harmonic {A}nalysis and {B}oundary {V}alue {P}roblems in the
  {C}omplex {D}omain}.   Birkh\"auser Verlag, Basel (1993).

\bibitem{Driscoll2014}
T.A. Driscoll, N. Hale, L. Trefethen,
{\em Chebfun Guide}. Pafnuty Publications (2014).

\bibitem{druskin1998extended}
V. Druskin, L. Knizhnerman,
Extended Krylov subspaces: approximation of the matrix square root
and related functions. {\em SIAM J. Matrix Anal. Appl.},
\textbf{19}, No 3 (1998), 755--771, DOI: https://doi.org/10.1137/S0895479895292400.

\bibitem{DLP-Pade}
B. Duan, R.D. Lazarov, J.E. Pasciak,
Numerical approximation of fractional powers of elliptic operators.
{\em IMA J. Numerical Analysis} \textbf{40}, No 3 (2019), 1746--1771;
DOI: 10.1093/imanum/drz013.

\bibitem{DU2020numerical}
M. D'Elia, Q. Du, C. Glusa, M. Gunzburger, X. Tian, and Z. Zhou, 
Numerical methods for nonlocal and fractional models. {\em Preprint arXiv:2002.01401} (2020).

\bibitem{Farago-05}
I. Farag\'{o}, Splitting methods and their application to the abstract Cauchy problems.
In: {\em  Numerical Analysis and Its Applications. NAA 2004. Lecture Notes in Computer Science}
(Eds: Z. Li, L. Vulkov, J. Waśniewski) \textbf{3401} (2005), 33--45.

\bibitem{peridynamics}
W.~H. Gerstle,  {\em Introduction to Practical Peridynamics}.
World Scientific, 2015.

\bibitem{GHH-19}
I. Georgieva, S. Harizanov, C. Hofreither,
Iterative low-rank approximation solvers for the extension method for
fractional diffusion.
{\em Computers $\&$ Mathematics with Appl.} \textbf{80}, No 2 (2020), 351--366;
DOI: 10.1016/j.camwa.2019.07.016.

\bibitem{gilboa2008nonlocal}
G.~Gilboa and S.~Osher,  Nonlocal operators with applications to image processing.
{\em Multiscale Modeling \& Simulation}, \textbf{7}, No 3 (2008), 1005--1028, DOI: 10.1137/070698592.

\bibitem{grubb2016}
G.~Grubb, Regularity of spectral fractional Dirichlet and Neumann problems.
{\em Mathematische Nachrichten} \textbf{289}, No 7 (2016), 831--844.

\bibitem{HLMM-19}
S. Harizanov, R. Lazarov, S. Margenov, P. Marinov,
{The best uniform rational approximation BURA of $t^\alpha$,
$t \in [0,1]$, $\alpha \in (0,1)$: Applications to solving equations
involving fractional powers of elliptic operators}.
\emph{Lecture Notes in Computer Science and Technologies}, No 9, IICT-BAS (2019).

\bibitem{HLMMV18}
S. Harizanov, R. Lazarov, S. Margenov, P. Marinov, Y. Vutov,
Optimal solvers for linear systems with fractional powers of sparse  SPD matrices.
{\em Numerical Linear Algebra with Applications} \textbf {25}, No 4 (2018),
115--128, DOI: 10.1002/nla.2167.

\bibitem{harizanov2019cmwa}
S. Harizanov, R. Lazarov, S. Margenov, P. Marinov,
Numerical solution of fractional diffusion–reaction problems based on BURA.
{\em Computers $\&$ Mathematics with Appl.}, \textbf{80}, No 2, (2020), 316--331;
DOI: 10.1016/j.camwa.2019.07.002.

\bibitem{HLMMP-19}
S. Harizanov, R. Lazarov, S. Margenov, P. Marinov, J. Pasciak,
Comparison analysis of two numerical methods for fractional diffusion problems
based on the best rational approximations of $t^\gamma$ on $[0,1]$.
{\em Lecture Notes in Computational Science and Engineering} \textbf{128} (2019), 165--185.

\bibitem{harizanov2020JCP}
S. Harizanov, R. Lazarov, S. Margenov, P. Marinov, J. Pasciak.
Analysis of numerical methods for spectral fractional elliptic
equations based on the best uniform rational approximation.
\newblock {\em J. of Computational Physics} \textbf{408} (2020);
DOI: 10.1016/j.jcp.2020.109285.

\bibitem{HatanoHatano:1998}
Y. Hatano and N. Hatano, Dispersive transport of ions in column experiments: an explanation of
  long-tailed profiles. {\em Water Resources Res.} \textbf{34} (1998), 1027--1033.

\bibitem{higham2008functions}
N.J. Higham,
{\em Functions of Matrices: Theory and Computation}. SIAM, 2008.

\bibitem{HOFREITHER2019}
C. Hofreither,
A unified view of some numerical methods for fractional diffusion.
{\em Computers \& Mathematics with Applications} \textbf{80}, No 2 (2020), 332--350.
DOI: 10.1016/j.camwa.2019.07.025.

\bibitem{HOFREITHER2020}
C. Hofreither,
An algorithm for best rational approximation based on barycentric rational interpolation.
{\em RICAM-Report} No 2020 -- 37 (2020).

\bibitem{ilic2009numerical}
M. Ili\'{c}, I.W. Turner, V. Anh, A numerical solution using an adaptively
preconditioned Lanczos method for a class of linear systems related with
the fractional Poisson equation. {\em Int. J. Stochastic Analysis}
\textbf{2008} (2009); DOI:10.1155/2008/104525.

\bibitem{kato}
T. Kato,
Fractional powers of dissipative operators. {\em J. Math. Soc. Japan}
\textbf{13}, No 3 (1961), 246--274.

\bibitem{KMV-18}
N. Kosturski, S. Margenov, Y. Vutov,
Performance Analysis of MG Preconditioning on Intel Xeon Phi:
Towards Scalability for Extreme Scale Problems with Fractional Laplacians.
In: {\em  Large-Scale Sci. Computing. LSSC 2017. Lecture Notes in Computer Sci., Springer, Cham} 
(Eds: I. Lirkov, S. Margenov) \textbf{10665} (2018), 304--312.

\bibitem{kwasnicki2017ten}
M.~Kwa{\'s}nicki, Ten equivalent definitions of the fractional Laplace operator.
{\em Fract. Calc. Appl. Anal.} \textbf{20}, No 1, (2017), 7--51; DOI: 10.1515/fca-2017-0002.

\bibitem{LPGSGZMCMAK-18}
A.~Lischke, G.~Pang, M.~Gulian, F.~Song, C.~Glusa,
X.~Zheng, Z.~Mao, W.~Cai, M.M.~Meerschaert, M.~Ainsworth, G.~Karniadakis,
What is the fractional Laplacian? A comparative review with new results.
{\em J. of Computational Physics} \textbf{404} (2020); DOI: 10.1016/j.jcp.2019.109009.

\bibitem{Marchuk-68}
G.I. Marchuk, Some applications of splitting-up methods to the solution
of problems in mathematical physics. {\em Aplikace Matematiky} \textbf{1} (1968), 103--132.

\bibitem{PGMASA1987}
P.G. Marinov, A.S. Andreev,
A modified Remez algorithm for approximate determination of the
rational function of the best approximation in Hausdorff metric.
{\em C.R. Acad. Bulg. Sci.} \textbf{40}, No 3 (1987), 13--16.

\bibitem{MRAplus-19}
S. Margenov,T. Rauber, E. Atanassov, F. Almeida, V. Blanco, R. Ciegis, A. Cabrera, N. Frasheri,
S. Harizanov, R. Kriauzien, G. Ruenger, P. San Segundo, V. Starikovicius, S. Szabo, B. Zavalnij,
Applications for ultra-scale systems.
{\em IET Professional Applications of Computing Series} \textbf{24} (2019), 189--244.

\bibitem{MetzlerKlafter:2004}
R.~Metzler and J.~Klafter, The restaurant at the end of the random walk: recent developments in
  the description of anomalous transport by fractional dynamics.  {\em J. Phys. A}, \textbf{37}, No 31 (2004), R161--R208.

\bibitem{NST-18}
Y. Nakatsukasa, O. S\'{e}te, L.N. Trefethen, The AAA algorithm for rational approximation,
{\em SIAM J. Sci. Comp.} \textbf{40}, No 3 (2018), A1494--A1522, DOI: 10.1137/16M1106122.

\bibitem{nochetto2015pde}
R.H. Nochetto, E. Ot{\'a}rola, A.J. Salgado,
A PDE approach to fractional diffusion in general domains: a priori error analysis.
{\em Found. Comput. Math.} \textbf{15}, No 3 (2015), 733--791.

\bibitem{nochetto2016pde}
R.H. Nochetto, E Ot{\'a}rola, A.J. Salgado,
A PDE approach to space-time fractional parabolic problems.
{\em SIAM J. Numer. Anal.} \textbf{54}, No 2 (2016), 848--873.

\bibitem{P-geophysical}
J. Pedlosky,
{\em  Geophysical Fluid Dynamics}. {Springer Science \& Business Media}, 2013.

\bibitem{Podlubny:1999}
I.~Podlubny,  {\em Fractional {D}ifferential {E}quations}.
{Acad. Press}, San Diego, CA, 1999.

\bibitem{Ros_Oton_2014}
X.~Ros-Oton and J.~Serra,  The {P}ohozaev identity for the fractional {L}aplacian.
{\em Archive for Rational Mechanics and Analysis} \textbf{213}, No 2, (2014), 587–-628.

\bibitem{saff1992asymptotic}
E.B. Saff, H. Stahl,
{\em Asymptotic Distribution of Poles and Zeros of Best Rational
Approximants to $x^\alpha$ on $[0,1]$}.
Topics in Complex Analysis, Banach Center Publ.,
Vol. 31, Institute of Mathematics, Polish Academy of Sciences, Warsaw (1995).

\bibitem{samarskii2001theory}
A.A. Samarskii,
{\em The Theory of Difference Schemes}. Ser. Pure and Applied Mathematics,
Vol. 240, Marcel Dekker, Inc., New York (2001).

\bibitem{SXK-17}
F. Song, C. Xu, G.E. Karniadakis,
Computing fractional laplacians on complex-geometry domains:
Algorithms and simulations.
{\em SIAM J. Sci. Comp.}, \textbf{39}, No 4 (2017), A1320--A1344.

\bibitem{Stahl93}
H. Stahl,
Best uniform rational approximation of $x^\alpha$ on $[0, 1]$.
{\em  Bull. Amer. Math. Soc. (N.S.)} \textbf{28}, No 1 (1993), 116--122.

\bibitem{Stahl2003}
H.R. Stahl,
Best uniform rational approximation of $x^\alpha$ on $[0, 1]$.
{\em Acta Math.} \textbf{190}, No 2 (2003), 241--306.

\bibitem{Strang-68}
G. Strang, On the construction and comparison of difference schemes.
{\em SIAM J. Num. Anal.} \textbf{5} (1968), 506-517.

\bibitem{SZBCC-18}
H.G. Sun, Y. Zhang, D. Baleanu, W. Chen, Y.Q. Chen,
A new collection of real world applications of fractional
calculus in science and engineering.
{\em Commun. Nonlin. Sci. Numer. Simul.} \textbf{64} (2018), 213--231.

\bibitem{mathworks09}
The MathWorks, Numerics::fMatrix--functional calculus for numerical square
matrices;
%\url
\href{http://www.mathworks.com/access/helpdesk/help/toolbox/mupad/numeric/fMatrix.html}
{http://www.mathworks.com/access/helpdesk/help/toolbox/mupad/numeric/%} \hfill \break \hspace*{0.1cm } \hfill {
fMatrix.html},  (2009).

\bibitem{Thomee2006}
V. Thom\'ee,
{\em Galerkin Finite Element Methods for Parabolic Problems}.
Springer Ser. in Comput. Mathematics, Vol. 25,
Springer-Verlag, Berlin, 2nd Ed. (2006).

\bibitem{Tref-20}
L.N. Trefethen, Y. Nakatsukasa, J.A.C. Weideman,
Exponential node clustering at singularities for rational
approximation, quadrature, and PDEs.
{\em arXiv:2007.11828v1} (2020).

\bibitem{Vabishchevich14}
P.N. Vabishchevich,
Numerical solving the boundary value problem for fractional
powers of elliptic operators.
{\em CoRR}, \textbf{abs/1402.1636} (2014).

\bibitem{Vabishchevich15}
P.N. Vabishchevich,
Numerically solving an equation for fractional powers of elliptic operators.
{\em J. of Comput. Phys.} \textbf{282} (2015), 289--302.

\bibitem{V-16}
P.N. Vabishchevich,
Numerical solution of non-stationary problems for a space-fractional
diffusion equation.
{\em Fract. Calc. Appl. Anal.} \textbf{19}, No 1 (2016), 116--139; 
DOI: 10.1515/fca-2016-0007.

\bibitem{Vabishchevich18}
P.N. Vabishchevich,
Numerical solution of time-dependent problems with fractional power
elliptic operator.
{\em Comput. Methods in Appl. Math.} \textbf{18}, No 1 (2018), 111--128.

\bibitem{V-20}
P.N. Vabishchevich,
Approximation of a fractional power of an elliptic operator.
{\em Numer. Lin.  Algebra with Appl.} \textbf{27}, No 3 (2020);
DOI: 10.1002/nla.2287.

\bibitem{varga1992some}
R.S. Varga, A.J. Carpenter,
Some numerical results on best uniform rational approximation of
$x^\alpha$ on [0, 1].
{\em Numerical Algorithms} \textbf{2}, No 2 (1992), 171--185.

\bibitem{Yanenko-62}
N.N. Yanenko, On convergence of the splitting method for heat equation with
variable coefficients. {\em J. Comput. Math. Math. Phys.}
\textbf{2}, No 5 (1962), 933--937 (in Russian).

\bibitem{XuZ2002}
J. Xu, L. Zikatanov,
The method of alternating projections and the method of subspace
  corrections in {H}ilbert space.
  {\em J. Amer. Math. Soc.} \textbf{15}, No 3 (2002), 573--597, 
  DOI: https://doi.org/10.1090/S0894-0347-02-00398-3.

\end{thebibliography}
\end{document}